\documentclass[a4paper,12pt]{amsart}
\usepackage{amsmath}
\usepackage{amsthm}
\usepackage{amscd}
\usepackage{amssymb}
\usepackage{amsfonts}
\usepackage{latexsym}
\usepackage[mathscr]{eucal}

\newtheorem{thm}{Theorem}[section]
\newtheorem{prop}[thm]{Proposition}
\newtheorem{lemma}[thm]{Lemma}
\newtheorem{cor}[thm]{Corollary}

\theoremstyle{definition}
\newtheorem{defi}[thm]{Definition}
\newtheorem{exam}[thm]{Example}
\theoremstyle{remark}
\newtheorem{rema}[thm]{Remark}
\numberwithin{equation}{section}

\def\graph{\mathop{\rm {graph}}\nolimits}

\def\graph{\mathop{\rm {graph}}\nolimits}

\def\Aut{\hbox{\rm Aut}}

\newcommand{\C}{\mathbb{C}}

\newcommand{\N}{\mathbb{N}}

\newcommand{\R}{\mathbb{R}}

\newcommand{\Z}{\mathbb{Z}}

\newcommand{\cB}{{\mathcal B}}
\newcommand{\cC}{{\mathcal C}}
\newcommand{\cT}{{\mathcal T}}

\newcommand{\cF}{{\mathcal F}}
\newcommand{\cG}{{\mathcal G}}
\newcommand{\cO}{{\mathcal O}}

\newcommand{\tf}{\tilde{f}}
\newcommand{\ta}{\tilde{a}}

\newcommand{\tS}{\tilde{S}}
\newcommand\inpr[2]{\langle{#1,#2}\rangle}
\newcommand\ORC{{\cO}_R(\hat{\C})} 
\begin{document}
\title{KMS states and branched points}
\author{Masaki Izumi}
\address[Masaki Izumi]{Department of Mathematics, 
Graduate School of Sciences, Kyoto University,
Kyoto 606-8502, Japan}

\author{Tsuyoshi Kajiwara}
\address[Tsuyoshi Kajiwara]{Department of Environmental and 
Mathematical Sciences, 
Okayama University, Tsushima, 700-8530,  Japan}      

\author{Yasuo Watatani}
\address[Yasuo Watatani]{Department of Mathematical Sciences, 
Kyushu University,
Hakozaki, Fukuoka, 812-8581,  Japan}

\maketitle
\begin{abstract}
We completely classify the KMS states for the gauge action on a 
$C^*$-algebra associated with a rational function $R$ introduced 
in our previous work. 
The gauge action has a phase transition at $\beta = \log \deg R$. 
We can recover the degree of $R$, the number of branched points, 
the number of exceptional points and the orbits of exceptional points 
from the structure  of the KMS states. 
We also classify the KMS states for $C^*$-algebras associated with 
some self-similar sets, including the full tent map 
and the Sierpinski gasket by a similar method. 

\end{abstract}

\section{Introduction}
It is of fundamental interest in operator algebras to analyze interplay 
between a geometric or dynamical object and a $C^*$-algebra associated 
with it. 
For a branched covering, Deaconu and Muhly \cite{DM} introduced a 
$C^*$-algebra associated with it using a r-discrete groupoid. 
A typical example of a branched covering is a rational function regarded 
as a self-map of the Riemann sphere $\hat{\mathbb C}$. 
In order to capture information of the branched points for the complex 
dynamical system arising from a rational function $R$,  
the second and third-named authors \cite{KW1} introduced a 
slightly different construction of a $C^*$-algebra 
${\mathcal O}_R(\hat{\mathbb C})$ 
(resp. ${\mathcal O}_R(J_R)$ and ${\mathcal O}_R(F_R)$) 
associated with $R$ on $\hat{\mathbb C}$ 
(resp. the Julia set $J_R$ and the Fatou set $F_R$ of $R$).    
The $C^*$-algebra ${\mathcal O}_R(\hat{\mathbb C})$ is the Cuntz-Pimsner 
algebra of a Hilbert bimodule over the $C^*$-algebra $C(\hat{\C})$ 
of the set of continuous functions and the other two are defined in a 
similar way. 

One of the purposes of the present paper is to discuss KMS states for 
the gauge action on the $C^*$-algebra $\ORC$. 
The structure of the KMS states reflects that of the singular points of $R$ 
as we expected.  
We completely classify the KMS states for the gauge action of 
${\mathcal O}_R(\hat{\mathbb C})$. 
If $R$ has no exceptional points, then the gauge action has a phase transition 
at $\beta = \log \deg R$ in the following sense:  
In the region $0 \leq \beta < \log \deg R$, no KMS-state exists.  
A unique KMS-state exists at $\beta = \log \deg R$, which is of type 
$III_{1/\deg R}$ and corresponds to the Lyubich measure.  
The extreme $\beta$-KMS states at $\beta > \log \deg R$ 
are parameterized by the branched points of $R$ and are factor states of 
type I.
If $R$ has exceptional points, then there appear additional $\beta$-KMS 
states for $0<\beta \leq \log \deg R$ parameterized by exceptional points. 
We can recover the degree of $R$, the number of branched points, 
the number of exceptional points  from the structure  of the KMS states. 
The orbits of exceptional points are distinguished by 0-KMS states. 

We also classify the KMS states for the $C^*$-algebras associated with some 
self-similar sets including the  full tent map and the Sierpinski gasket by 
a similar method. 

Olsen-Pedersen \cite{OP} showed that a $\beta$-KMS state for the gauge action 
of the Cuntz algebra ${\mathcal O}_n$ exists if and only if $\beta = \log n$ 
and that $\log n$-KMS state is unique. 
Since then, several authors have discussed KMS states for the gauge action 
(and its generalization) of the Cuntz-Pimsner algebra 
(and its generalization). 
Here we content ourselves with only giving an incomplete list of 
such works: 
\cite{DS}, \cite{EFW}, \cite{Ev}, \cite{Ex1}, \cite{Ex2}, \cite{EL}, \cite{I}, 
\cite{KP}, \cite{KR}, \cite{LN}, \cite{M}, \cite{MWY}, \cite{OK}, \cite{PWY}. 
Among the others, in this paper we follow Laca and Neshveyev's approach where 
the structure of the KMS states is described in terms of a certain 
Perron-Frobenius type operator. 
We give an explicit description of the Perron-Frobenius type operator 
in our cases, which allows us to perform detailed analysis of the KMS states. 

This paper is an extended version of the preprint \cite{KW3}.

\section{Dynamical systems and Hilbert bimodules}

In this section, we recall our construction of the $C^*$-algebras associated 
with rational functions in \cite{KW1} and self-similar sets in \cite{KW2}, 
which are constructed as the Cuntz-Pimsner algebras \cite{Pi}.    

\subsection{The Cuntz-Pimsner algebra}
Let $A$ be a $C^*$-algebra and $X$ be a Hilbert right $A$-module.  
We denote by $L(X)$ the algebra of the adjointable bounded operators 
on $X$.  
For $\xi$, $\eta \in X$, the ``rank one" operator $\theta _{\xi,\eta}$
is defined by $\theta _{\xi,\eta}(\zeta) = \xi(\eta|\zeta)_A$
for $\zeta \in X$. 
The closure of the linear span of the rank one operators is denoted by $K(X)$. 

A sequence $(u_n)_n$ in $X$ is said to be a countable basis of $X$ over $A$ 
if for any $x \in X$, the series 
$\sum_{n=1}^{\infty}  u_n(u_n|x)_A$ converges to $x$ in norm. 
We note that $(u_n)_n$ converges unconditionally in the 
sense that the net  $(\sum_{n \in  F}   u_n(u_n|x)_A)_F$, 
where $F$ runs  all finite subsets of ${\mathbb N}$, converges 
to $x$ in norm. (See \cite{KPW} on basis). 
We have  $\|u_n \| \leq 1$ and the sequence 
$(\sum_{k=1}^n \theta _{u_k,u_k})_n$ is an approximate unit for $K(X)$.     

We say that 
$X$ is a Hilbert bimodule over $A$ if $X$ is a Hilbert right  $A$-
module with a *-homomorphism $\phi : A \rightarrow L(X)$.  We always assume 
that $X$ is full and $\phi$ is injective. 
   Let $F(X) = \bigoplus _{n=0}^{\infty} X^{\otimes n}$
be the full Fock module of $X$ with a convention $X^{\otimes 0} = A$. 
 For $\xi \in X$, the creation operator $T_{\xi} \in L(F(X))$ is defined by 
\[
T_{\xi}(a) =  \xi a  \qquad \text{and } \ 
T_{\xi}(\xi _1 \otimes \dots \otimes \xi _n) = \xi \otimes 
\xi _1 \otimes \dots \otimes \xi _n .
\]
We define $i_{F(X)}: A \rightarrow L(F(X))$ by 
$$
i_{F(X)}(a)(b) = ab \qquad \text{and } \ 
i_{F(X)}(a)(\xi _1 \otimes \dots \otimes \xi _n) = \phi (a)
\xi _1 \otimes \dots \otimes \xi _n 
$$
for $a,b \in A$.  The Cuntz-Toeplitz algebra ${\mathcal T}_X$ 
is the C${}^*$-algebra acting on $F(X)$ generated by $i_{F(X)}(a)$
with $a \in A$ and $T_{\xi}$ with $\xi \in X$. 

Let $j_K | K(X) \rightarrow {\mathcal T}_X$ be the homomorphism 
defined by $j_K(\theta _{\xi,\eta}) = T_{\xi}T_{\eta}^*$. 
We consider the ideal $I_X := \phi ^{-1}(K(X))$ of $A$. 
Let ${\mathcal J}_X$ be the ideal of ${\mathcal T}_X$ generated 
by $\{ i_{F(X)}(a) - (j_K \circ \phi)(a) ; a \in I_X\}$.  Then 
the Cuntz-Pimsner algebra ${\mathcal O}_X$ is defined as 
the quotient ${\mathcal T}_X/{\mathcal J}_X$ . 
Let $\pi : {\mathcal T}_X \rightarrow {\mathcal O}_X$ be the 
quotient map.  
We set $S_{\xi} = \pi (T_{\xi})$ and $i(a) = \pi (i_{F(X)}(a))$. 
Let $i_K : K(X) \rightarrow {\mathcal O}_X$ be the homomorphism 
defined by $i_K(\theta _{\xi,\eta}) = S_{\xi}S_{\eta}^*$. Then 
$\pi((j_K \circ \phi)(a)) = (i_K \circ \phi)(a)$ for $a \in I_X$.   

The Cuntz-Pimsner algebra ${\mathcal O}_X$ is 
the universal C${}^*$-algebra generated by $i(a)$ with $a \in A$ and 
$S_{\xi}$ with $\xi \in X$  satisfying that 
$i(a)S_{\xi} = S_{\phi (a)\xi}$, $S_{\xi}i(a) = S_{\xi a}$, 
$S_{\xi}^*S_{\eta} = i((\xi | \eta)_A)$ for $a \in A$, 
$\xi, \eta \in X$ and $i(a) = (i_K \circ \phi)(a)$ for $a \in I_X$.
We usually identify $i(a)$ with $a$ in $A$.  
We also identify $S_{\xi}$ with $\xi \in X$ and simply write $\xi$
instead of $S_{\xi}$.  

There exists an action 
$\alpha : {\mathbb R} \rightarrow \Aut \ {\mathcal O}_X$
defined by $\alpha_t(\xi) = e^{it}\xi$ for $\xi\in X$ and $\alpha_t(a)=a$ 
for $a\in A$, which is called the {\it gauge action}.

%
%
%
 
\subsection{The case of rational functions}  
Let $R$ be a rational function with $N = \deg R$, that is,   
if $R(z)=P(z)/Q(z)$ with relatively prime polynomials $P(z)$ and $Q(z)$, 
the degree $\deg R$ is the maximum of those of $P(z)$ and $Q(z)$.  
We regard $R$ as a $N$-fold branched covering map  
$R : \hat{\mathbb C} \rightarrow \hat{\mathbb C}$
on the Riemann sphere $\hat{\mathbb C} = {\mathbb C} 
\cup \{ \infty \}$.  
The sequence $(R^n)_n$ of iterations of $R$ 
gives a complex dynamical system on $\hat{\mathbb C}$. 
The Fatou set $F_R$ of $R$ is the maximal open subset of 
$\hat{\mathbb C}$ on which $(R^n)_n$ is equicontinuous (or 
a normal family), and the Julia set $J_R$ of $R$ is the 
complement of the Fatou set in $\hat{\mathbb C}$.  
We always assume that $R$ is of degree at least two.  

Recall that a {\it branched point} of $R$ is a point   
$z_0$ around which $R$ is not locally one to one. 
It is a zero of the derivative $R'$ or 
a pole of $R$ of order two or higher.  
The image $w_0 =R(z_0)$ is called a {\it branch value} of $R$.    
Using appropriate local charts, if $R(z) = w_0 + c(z - z_0)^n + 
(\text{higher terms})$ with 
$n \geq 1$ and $c \not= 0$ on some neighborhood of $z_0$, 
then the integer $n = e(z_0)$ is called the 
{\it branch index} of $R$ at $z_0$.  
Thus $e(z_0) \geq 2$ if $z_0$ is a branched point and $e(z_0) = 1$ 
otherwise. 
Therefore $R$ is an $e(z_0) :1$ map in a punctured neighborhood of $z_0$.
Let $\cB(R)$ be the 
set of branched points of $R$ and $\cC(R)$ be the 
set of the branch values of $R$. 
Then the restriction $R: \hat{\mathbb C} \setminus R^{-1}(\cC(R))
\rightarrow \hat{\mathbb C} \setminus \cC(R) $ is an 
$N$-to-one regular covering. 

Let $A= C(\hat{\mathbb C})$ and $X = C(\graph R)$ be the set of continuous 
functions on $\hat{\mathbb C}$ and $\graph R$ respectively, 
where $\graph R = \{(x,y) \in \hat{\mathbb C}^2 ; y = R(x)\} $ 
is the graph of $R$.
Then $X$ is an $A$-$A$ bimodule by 
$$
(a\cdot \xi \cdot b)(x,y) = a(x)\xi(x,y)b(y),\quad a,b \in A,\; 
\xi \in X.$$ 
We define an $A$-valued inner product $(\ |\ )_A$ on $X$ by 
$$
(\xi|\eta)_A(y) = \sum _{x \in R^{-1}(y)} e(x) \overline{\xi(x,y)}\eta(x,y),
\quad \xi,\eta \in X,\; y \in \hat{\mathbb C}.$$  
Thanks to the branch index $e(x)$, the inner product above gives a continuous  
function and $X$ is a full Hilbert bimodule over $A$ without completion. 
The left action of $A$ is unital and faithful. 

Since the Julia set $J_R$ is completely invariant under $R$, i.e., 
$R(J_R) = J_R = R^{-1}(J_R)$, we can consider the restriction 
$R|_{J_R} : J_R \rightarrow J_R$, which will be often denoted by 
the same letter $R$.
 Let $\graph R|_{J_R} = \{(x,y) \in J_R \times J_R \ ; \ y = R(x)\} $
be the graph of the restriction map $R|_{J_R}$ and 
$X(J_R) = C(\graph R|_{J_R})$. 
In the same way as above, $X(J_R)$ is a full Hilbert bimodule over $C(J_R)$. 
Since the Fatou set $F_R$ is also completely invariant, 
$X(F_R):=C_0(\graph R|_{F_R})$ is a full Hilbert bimodule over $C_0(F_R)$.

\begin{defi}[\cite{KW1}]   
The $C^*$-algebra
${\mathcal O}_R(\hat{\mathbb C})$ is defined 
as the Cuntz-Pimsner algebra of the Hilbert bimodule 
$X= C(\graph R)$ over 
$A = C(\hat{\mathbb C})$. 
When the Julia set $J_R$ is not empty (for example 
$\deg R \geq 2$), we define 
the $C^*$-algebra ${\mathcal O}_R(J_R)$ 
as the Cuntz-Pimsner algebra of the Hilbert bimodule $
X= C(\graph R|_{J_R})$ over $A = C(J_R)$. 
When the Fatou set $F_R$ is not empty, the $C^*$-algebra
${\mathcal O}_R(F_R)$ is defined similarly.
\end{defi}

\subsection{The case of self-similar sets} 
Let $(\Omega,d)$ be a separable, complete metric space $\Omega$ with a 
metric $d$. Let $\gamma = (\gamma_1,\dots , \gamma_N)$ be a system of    
continuous maps from $\Omega$ to itself. We say that $\gamma_i$ is 
a proper contraction if  
there exist positive constants $c_1$ and $c_2$
with $0 < c_1 \leq c_2 < 1$ satisfying the condition:
$$c_1 d(x,y)\leq d(\gamma_i(x),\gamma_i(y))\leq c_2 d(x,y), 
\quad  i=1,2,..,N,\quad \forall x,y \in X.$$ 

We say that a non-empty compact set $K \subset \Omega$ is {\it self-similar}
(in a weak sense) with respect to the system 
$\gamma = (\gamma_1,\dots , \gamma_N)$ if 
$$
K=\bigcup_{i=1}^{N} \gamma_{i}(K).
$$
If the  contractions are  proper, then there exists a unique self-similar 
set $K \subset X$. 
In this note  we usually forget the ambient space $\Omega$ and 
assume that $\gamma = (\gamma_1,\dots , \gamma_N)$ is a 
system of continuous functions on a self-similar $K$.

We use the following notations:  
\begin{align*} 
 \cB(\gamma) & = \{ x \in K | x =\gamma_j(y)=\gamma_{j'}(y) 
\text{ for some } y 
 \in K \text{ and } j \ne j' \} \\
 \cC(\gamma) & =  \{ y \in K | \gamma_j(y)=\gamma_{j'}(y) \text{ for some } 
 j \ne j' \} 
\end{align*}
We call a point in $\cB(\gamma)$ a branched point,  and a point in
$\cC(\gamma)$ a branch value. 
If  $\gamma = (\gamma_1,\dots , \gamma_N)$ is a system 
of branches of $R^{-1}$ for a certain map $R$, the terms are   
compatible with those for $R$.  


We set $\cG$ to be the union of the cographs of 
$\gamma_i$ for $i=1,2,\cdots,N$, that is, 
\[
{\mathcal G} = {\mathcal G}(\{\gamma_j:j=1,2,..,N \}) := 
\bigcup _{i=1}^N \{(x,y) \in K^2 ; x = \gamma _i(y)\}. 
\] 
Let $A = C(K)$ and let $X = C({\mathcal G})$.  
Then $X$ is an $A$-$A$ bimodule by 
$$
(a\cdot f \cdot b)(x,y) = a(x)f(x,y)b(y),\quad 
a,b \in A,\;f \in X.$$ 
We introduce an $A$-valued inner product $(\ |\ )_A$ on $X$ by 
$$
(\xi|\eta)_A(y) = \sum _{i=1}^N  
\overline{\xi(\gamma _i(y),y)}\eta(\gamma _i(y),y),\quad 
\xi,g \in X,\; y \in K.$$ 

\begin{defi}[\cite{KW2}]  
Let $(K,d)$ be a compact metric space 
and $\gamma = (\gamma_1,\dots , \gamma_N)$ be a system of 
proper contractions on $K$.  Assume that $K$ is self-similar. 
The $C^*$-algebra ${\mathcal O}_{\gamma}(K)$ is defined 
as the Cuntz-Pimsner algebra ${\mathcal O}_X$ of the Hilbert bimodule 
$X= C({\mathcal G})$ over $A = C(K)$. 
\end{defi}

\section{KMS states}
First we recall some facts on KMS states on general 
Cuntz-Pimsner algebras from Laca and Neshveyev \cite{LN}. 
Let $A$ be a $C^*$-algebra and $X$ be a full Hilbert 
$A$-$A$ bimodule with a non-degenerate left $A$-action. 
Let $\alpha$ be the gauge action. 
For $\beta>0$, we denote by $K_\beta(\alpha)$ the set of $\beta$-KMS states 
for $(\cO_X,\alpha)$.

In the following we assume that there exists a countable basis 
$\{u_i\}_{i=1}^\infty$ for $X$ as a right Hilbert $A$-module 
and assume that  
$N:=\sup_{n\in \N}||\sum_{i=1}^n(u_i|u_i)_A||$ is finite.

\begin{defi}[Perron-Frobenius type operators] 
Let the notation be as above. 
We introduce a Perron-Frobenius type operator $F:A^*\rightarrow A^*$, 
which is a bounded positive map, by the following formula: 
$$
F(\omega)(a)= \sum_{i=1}^\infty \omega((u_i|a \cdot u_i)_A),
\quad \omega\in A^*,\;a\in A.$$ 
For $\beta\geq 0$, we set $F_\beta = e^{-\beta} F$.  
\end{defi}

The operator $F$ was discussed in \cite{PWY} in the case where $X$ is finitely generated  
and in \cite{LN} in the general case. 
When $\tau$ is  a finite positive trace on $A$,  it is known that $F(\tau)$ is again 
a finite positive trace, which does not depend on the choice of 
$\{u_i\}_{i=1}^\infty$. 
Therefore for positive $a\in A^+$ and a finite positive trace $\tau$, we have  
$$F(\tau)(a)=\sup_{v_i}\sum_{i=0}^n\tau((v_i|a\cdot v_i)),$$
where the supremum is taken over all $v_0,v_1, \cdots,v_n\in X$ satisfying 
$\sum_{i=1}^n\theta_{v_i,v_i}\leq I$ (see \cite[Theorem 1.1]{LN}).

\begin{thm}[Laca-Neshveyev \cite{LN}] 
Let $A$ and $X$ be as above and let $\alpha$ be the gauge action 
on $\mathcal O_X$.  
Then there is a bijective affine isomorphism
between $K_\beta(\alpha)$ and the set $T(A)_{\beta}$ of tracial states 
$\tau$ on $A$ satisfying 
the following conditions: 
$$F_\beta(\tau)(a)=\tau(a),\quad \forall a\in I_X,\leqno (K1)$$
$$F_\beta(\tau)(a)\leq \tau(a),\quad \forall a\in A^+.\leqno (K2)$$
The correspondence is given by restriction. 
\label{thm;K1K2}
\end{thm}

\begin{defi} Let $A$ and $X$ be as above. 
Following Exel and Laca \cite{EL}, we say that  
a finite positive trace $\tau _1$ on $A$ is of {\it finite type}  
if there exists a finite positive trace $\tau _0$ on $A$ such that 
$\tau _1 = \sum _{n=0}^{\infty} F_\beta^n(\tau_0)$ in the weak* 
topology. 
We denote by $T_f(A,F_\beta)$ the set of finite positive traces on 
$A$  of finite type. 
We say that a finite positive trace $\tau _2$ on $A$ is of 
{\it infinite type} if $F_\beta (\tau_2) = \tau _2$.  
We denote by $T_{i}(A,F_\beta)$ 
the set of finite positive traces on $A$  of infinite type.   

A $\beta$-KMS states $\varphi$ for the gauge action $\alpha$ on $\mathcal O_X$
is said to be of {\it finite type} (resp. {\it infinite type}) 
if the restriction $\tau=\varphi|_A$ is of finite type 
(resp. infinite type). 
We denote by $K_\beta(\alpha)_f$ (resp. $K_\beta(\alpha)_i$) 
the set of $\beta$-KMS states of finite (resp. infinite) type. 
\end{defi}

Laca and Neshveyev \cite[Proposition 2.4]{LN}  showed that 
any tracial state $\tau \in T(A)_{\beta}$ is uniquely decomposed as  
\begin{equation}
\tau= \tau_1 + \tau_2
\end{equation}
with $\tau_1 \in T_f(A,F_\beta)$ and 
$\tau_2 \in T_{i}(A,F_\beta)$.  
Moreover $\tau_1$ is given by 
\begin{equation}\tau_1 = \sum _{n=0}^{\infty} F_\beta^n(\tau_0)
\end{equation} 
with 
\begin{equation}\tau_0 = \tau - F_\beta(\tau)\end{equation} 
and $\tau_2$ is given by
\begin{equation}
\tau_2 = \lim _n F_\beta^n(\tau), \text{ in the weak* topology }. 
\end{equation}

The above classification of the KMS states a priori depends on the construction 
of $\cO_X$ from $A$ and $X$, and we don't know whether there is an intrinsic 
characterization of it in terms of the system $({\mathcal O}_X, \alpha)$.

Since the decomposition (3.1) is unique, we have 
\begin{equation}
\mathrm{ex}(K_\beta(\alpha))=\mathrm{ex}(K_\beta(\alpha)_f)\cup 
\mathrm{ex}(K_\beta(\alpha)_i).
\end{equation}
where for a convex set $C$, we denote by $\mathrm{ex}(C)$ 
the set of extreme points of $C$. 

(3.3) shows that $\tau_0$ vanishes on $I_X$ thanks to the condition (K1). 
Let $T(A/I_X)$ be the set of tracial states on $A$ that vanish  
on $I_X$, which may be regarded as tracial states on $A/I_X$.   
We denote by $T(A/I_X)_{\beta}$ the set of $\tau  \in T(A/I_X)$ 
such that 
$\sum_{n=0}^\infty F_\beta^n(\tau )$ 
converges in the weak$^*$-topology. 
When  $\beta>\log N$, we have $||F_\beta||<1$ and so  
$T(A/I_X)=T(A/I_X)_{\beta}$. 
For $\tau\in T(A/I_X)_{\beta}$, we set $\varphi_{\tau,\beta}$ to be 
the $\beta$-KMS state corresponding to the tracial state 
$\psi_{\tau,\beta}
:= m_{\tau,\beta}\sum_{n=0}^\infty F_\beta^n(\tau) \in T(A)_\beta,$
where $m_{\tau,\beta}$ is a normalizing constant. That is, 
$m_{\tau,\beta} = (\sum_{n=0}^\infty F_\beta^n(\tau)(1))^{-1}$ and 
$$
\varphi_{\tau,\beta}|_A 
= m_{\tau,\beta}\sum_{n=0}^\infty F_\beta^n(\tau) 
= \psi_{\tau,\beta}.
$$

\begin{lemma} The map 
$T(A/I_X)_{\beta}\ni \tau\mapsto \varphi_{\tau,\beta}\in K_\beta(\alpha)_f$ 
is a bijection, which is not necessarily affine,  but sends 
$\mathrm{ex}(T(A/I_X)_{\beta}) =
T(A/I_X)_{\beta}\cap \mathrm{ex}(T(A/I_X))$ onto $\mathrm{ex}(K_\beta(\alpha)_f)$. 
\label{lemma;finite}
\end{lemma}

\begin{proof} If a positive functional $\omega$ is 
dominated by a positive functional $\tau$ such that 
$\sum_{n=0}^\infty F_\beta^n(\tau )$ 
converges, then $\sum_{n=0}^\infty F_\beta^n(\omega )$ also
converges. Hence $\mathrm{ex}(T(A/I_X)_{\beta}) =
T(A/I_X)_{\beta}\cap \mathrm{ex}(T(A/I_X))$.  

The fact that the above map is bijective follows from the fact 
that $\tau$ is uniquely given by 
\[
\tau 
=((\varphi_{\tau,\beta}|_A - F_{\beta}(\varphi_{\tau,\beta}|_A)(1))^{-1}
(\varphi_{\tau,\beta}|_A - F_{\beta}(\varphi_{\tau,\beta}|_A)).
\] 
Let $\tau\in T(A/I_X)_\beta$ and $\tau', \tau''\in T(A/I_X)$ 
satisfy $\tau=c_1\tau'+ c_2\tau''$  
for some non-negative constants $c_1$ and $c_2$ with $c_1 + c_2 =1$ . 
Then $\tau'$ and $\tau''$ belong to $T(A/I_X)_\beta$, and we get 
$$
\varphi_{\tau,\beta}
=\frac{m_{\tau,\beta}c_1}{m_{\tau',\beta}}\varphi_{\tau',\beta}
+ \frac{m_{\tau,\beta}c_2}{m_{\tau'',\beta}}\varphi_{\tau'',\beta}.
$$  
Therefore if $\varphi_{\tau,\beta}$ is extreme, then so is $\tau$. 
Conversely assume that 
$\tau \in \mathrm{ex}(T(A/I_X)_{\beta})$ and that 
there exist $\varphi',\varphi''\in K_\beta(\alpha)_f$ such that 
$$\varphi_{\tau,\beta}=t_1\varphi'+t_2\varphi'', $$ 
for some  non-negative constants $t_1,t_2$ with 
$t_1 + t_2 =1$.  
Then there exists 
$\tau'$ and $\tau''\in T(A/I_X)_\beta$ 
such that $\varphi_{\tau',\beta}=\varphi'$ and 
$\varphi_{\tau'',\beta}=\varphi''$. Put 
$$c_1 = \frac{m_{\tau',\beta}t_1}{m_{\tau,\beta}}
\text{ and  }            
c_2 = \frac{m_{\tau'',\beta}t_2}{m_{\tau,\beta}}.
$$
Then 
\begin{align*}
& m_{\tau,\beta}\sum_n F_{\beta}^n(\tau)
= \varphi_{\tau,\beta}|A 
= t_1\varphi'|A+t_2\varphi''|A \\
&= m_{\tau,\beta}c_1\sum_n F_{\beta}^n(\tau') + 
  m_{\tau,\beta}c_2\sum_n F_{\beta}^n(\tau'') 
= m_{\tau,\beta}\sum_n F_{\beta}^n(c_1\tau' +c_2\tau''),
\end{align*}
which shows $\tau = c_1\tau' + c_2\tau''$.
Since $\tau$ is extreme, we get $\tau'=\tau''=\tau$ and so 
$\varphi'=\varphi''=\varphi_{\tau,\beta}$, which finishes the proof.
\end{proof}

\begin{cor} 
If $\beta>\log N$, then $\mathrm{ex}(K_\beta(\alpha))=
\{\varphi_{\tau,\beta}\}_{\tau\in \mathrm{ex}(T(A/I_X))}$. 
\label{cor;parametrization}
\end{cor}
\begin{proof}
If $\beta>\log N$, then $||F_\beta||<1$ holds, 
and so $K_\beta(\alpha)_i=\emptyset$. 
\end{proof}

In the rest of this section, we investigate the type of the GNS 
representation of a finite type $\beta$-KMS state $\varphi$ 
on $\cO_X$ for the gauge action. 

Let $(\pi_{\varphi},H_{\varphi},\Omega _{\varphi})$ 
be the  GNS representation associated with $\varphi$.  
Since $\{u_i\}_{i=1}^\infty$ is a countable basis for $X$, 
the sequence $(a_n)_n := (\sum_{k=1}^n \theta _{u_k,u_k})_n$ is an 
approximate unit for $K(X) \subset \cO_X$. 
Therefore $(\pi_{\varphi}(a_n))_n$ converges to a projection 
$p \in \pi_{\varphi}(\cO_X)''$. 
For $x\in X$, we have 
$$p\pi_\varphi(S_x)=\sum_{i=1}^\infty 
\pi_\varphi(S_{u_i(u_i|x)_A})=\pi_\varphi(S_x),  
\text{ in the strong topology } .$$
Thus for any $a\in A$
$$(1-p)\pi_{\varphi}(a)p=
\sum_{i=1}^\infty (1-p)\pi_{\varphi}(S_{au_i}S_{u_i}^*) =0,  
\text{ in the strong topology } .
$$
In particular, $p$ commutes with $\pi_{\varphi}(A)$.  
We denote by $\hat{\varphi}$ the weakly continuous extension of 
the state $\varphi$ to $\pi_{\varphi}(\cO_X)''$ given by 
$\hat{\varphi} (T) = (T\Omega_{\varphi} | \Omega_{\varphi})$ 
for $T \in \pi_{\varphi}(\cO_X)''$.  
 
\begin{lemma} Let $\varphi$ be a $\beta$-KMS state for $\alpha$ and 
let the notation be as above. 
Then the followings are equivalent:
\begin{itemize}
\item [$(1)$] $\varphi$ is of infinite type. 
\item [$(2)$] $\hat{\varphi} (p) = 1$. 
\end{itemize}
\label{lemma;p=1}
\end{lemma}
\begin{proof}Since $\varphi$ is a $\beta$-KMS state, 
for any $a \in A$, 
\begin{align*}
 (F_{\beta}(\varphi |_{A}))(a) 
& = \sum_{n=1}^\infty \frac{1}{e^{\beta}}\varphi |_{A}((u_n |au_n)_A)
  = \sum_{n=1}^\infty \frac{1}{e^{\beta}}\varphi (S_{u_n}^*aS_{u_n}) \\
& = \sum_{n=1}^\infty \varphi (S_{u_n}S_{u_n}^*a) 
  = \sum_{n=1}^\infty 
(\pi_{\varphi}(S_{u_n}S_{u_n}^*a)\Omega_{\varphi} | \Omega_{\varphi}) \\
& = (p\pi_{\varphi}(a)\Omega_{\varphi} | \Omega_{\varphi}) 
  = \hat{\varphi} (p\pi_{\varphi}(a)). 
\end{align*}
Putting $a = I$, we get $(F_{\beta}(\varphi |_{A}))(1) = \hat{\varphi} (p)$. 
Since $\varphi$ is a $\beta$-KMS state, the condition (K2) in Theorem \ref{thm;K1K2} 
implies $F_{\beta}(\varphi |_{A}) \leq  \varphi |_{A}$.  
Therefore $1 - \hat{\varphi} (p) = 0$ if and only if 
$\varphi |_{A}(1) - (F_{\beta}(\varphi |_{A}))(1) = 0$ 
if and only if $\varphi |_{A} = F_{\beta}(\varphi |_{A})$,  
i.e.,  $\varphi$ is of infinite type. 
\end{proof}

\begin{thm} In the above setting, 
let $\tau\in \mathrm{ex}(T(A/I_X)_{\beta})$ and 
$\varphi := \varphi_{\tau,\beta}$ be the corresponding extreme 
$\beta$-KMS state of finite type. 
Let $\tau'$ be a state on 
$A/I_X$ defined by $\tau'(a + I_X) = \tau (a)$.
Then $\pi_{\varphi}(\cO_X)''$ is of type I (resp. type II) if and only 
if $\pi_\tau'(A/I_X)"$ is of type I (resp. type II).  
In particular, if $A$ is abelian, then  $\varphi_{\tau,\beta}$ 
is a type I state. 
\end{thm}

\begin{proof} Since the  $\beta$-KMS state $\varphi$ is of 
finite type, $\hat{\varphi} (p) \not= 1$ by Lemma \ref{lemma;p=1}. 
Hence $p \not= I$. 
Let $H_0 = (I-p)H_{\varphi}$. 
Since $p$ commutes with $\pi_{\varphi}(A)$, 
we can  define a representation 
$\rho : A \rightarrow B(H_0)$ by 
$\rho (a) = (I-p)\pi_{\varphi}(a)|_{H_0}$. 
For any $x \in X^{\otimes n}$ 
and $y \in X^{\otimes m}$, we have $(I-p)\pi_{\varphi}(S_xS_y^*) = 0$. 
Therefore  
\[
(I-p)\pi_{\varphi}(\cO_X)''(I-p)
=(I-p)\pi_{\varphi}(A)''(I-p)=(I-p)\pi_{\varphi}(A)'' 
\cong  \rho (A)'' 
\]
Put $\Omega _0 := \frac{1}{\sqrt{1-\hat{\varphi} (p)}}(I-p)\Omega_{\varphi}$. 
Then $\|\Omega _0 \| = 1$.  For $a \in A$, 
\begin{align*}
& (\rho (a)\Omega _0|\Omega _0) 
= \frac{1}{1-\hat{\varphi} (p)}
((I-p)\pi_{\varphi}(a)\Omega_{\varphi}|\Omega_{\varphi}) 
= \frac{1}{1-\hat{\varphi} (p)}\hat{\varphi}((I-p)\pi_{\varphi}(a)) \\
& = \frac{1}{1-\hat{\varphi} (p)}
(\varphi |_{A}(a) - (F_{\beta}(\varphi |_{A}))(a)) 
= \tau (a). 
\end{align*}
Since $\rho (A)\Omega_0 = (I-p)\pi_{\varphi}(\cO_X)\Omega_{\varphi}$, 
the vector $\Omega_0$ is cyclic for $\rho$.  
Therefore $\rho$ is unitarily equivalent to the GNS representation 
$\pi _{\tau}$ of $A$ for $\tau$. 
Hence  $(I-p)\pi_{\varphi}(\cO_X)''(I-p)$ is isomorphic to 
$\pi _{\tau}(A)''$.  
Since $\pi_{\varphi}(\cO_X)''$ is a factor and $\pi _{\tau}(A)''$ is 
isomorphic to $\pi _{\tau '}(A/I_X)''$, we get the statement. 
\end{proof}

\section{Branched  points and Perron-Frobenius type operators}
In this section we give a detailed description of the Perron-Frobenius type 
operators $F$ introduced in the previous section for concrete examples. 
We treat the Hilbert bimodules arising from rational functions and 
self-similar sets in a unified setting, that is, a topological 
relation in the sense of Brenken \cite{Br} and a topological quiver with 
a weighted counting measure in the sense of Muhly-Solel \cite{MS} and 
Muhly-Tomforde \cite{MT}.  

Firstly we unify the notation of bimodules.  
Let $K$ be a compact metric space and $\cG$ be a closed subset of 
$K\times K$. 
For $(x,y)\in \cG$, we set $p_1(x,y)=x$ and $p_2(x,y)=y$, and 
$\cG_y=p_1p_2^{-1}(y)$. 
We assume $p_1(\cG)=p_2(\cG)=K$ and  
$N:=\sup_{y\in K} \#p_2^{-1}(y)<\infty.$  
Let $A=C(K)$, $X=C(\cG)$. 
We introduce an $A$-$A$ bimodule structure into $X$ by 
$(a\cdot \xi\cdot b)(x,y)=a(x)\xi(x,y)b(y).$ 
We assume that there exists a positive function $e$ on $\cG$ such 
that $e(x,y)\geq 1$ for all $(x,y)\in \cG$ and 
$$K\ni y\mapsto (\xi|\eta)_A(y):=\sum_{x\in \cG_y}e(x,y)
\overline{\xi(x,y)}\eta(x,y)$$
is continuous for all $\xi,\eta\in X$. 

\begin{lemma}\label{lemma;openneibourhood} 
For every $(x_0,y_0)\in \cG$, every open neighbourhood  
$W$ of $(x_0,y_0)$, and every $\varepsilon>0$, there exist open 
neighbourhoods $U$ of    
$x_0$ and $V$ of $y_0$ such that $U\times V\subset W$, 
$p_2((U\times V)\cap\cG)=V$, and 
$$|e(x_0,y_0)-\sum_{x\in U\cap \cG_y}e(x,y)|<\varepsilon,\quad  
\forall y\in V.$$ 
In particular, $e$ is upper semicontinuous. 
\end{lemma}

\begin{proof} Let $\cG_{y_0}=\{x_0,x_1,\cdots,x_n\}$. 
We choose open neighbourhoods $U_a$ of $x_a$ for $0\leq a\leq n$ and 
open neighbourhood $V_0$ of $y_0$ such that 
$U_0\times V_0\subset W$ and 
$\overline{U_a}\cap \overline{U_b}=\emptyset$ for $a\neq b$. 
Then there exists an open neighbourhood $V_1$ of $y_0$ 
with $ V_1\subset V_0$ such that 
for all $y\in V_1$,
$\cG_y\subset \bigcup_{a=0}^nU_a.$
Indeed, if there were no such $V_1$, there would exist a sequence 
$\{(x_n,y_n)\}_{n=1}^\infty$ in $\cG$ such that $\{y_n\}_{n=1}^\infty$ 
converges to $y_0$ and $x_n\notin U_a$ for all $a$. 
An accumulation point of this sequence is of the form $(x',y_0)\in \cG$ 
such that $x'\notin U_a$ for all $a$, which is a contradiction. 
Therefore $V_1$ as above exists. 
Let $\xi\in X$ such that $\xi(x,y)=1$ for $(x,y)\in U_0\times V_1$ and 
$\xi(x,y)=0,\quad \forall x\in \bigcup_{a=1}^nU_a.$
Then for $y\in V_1$, we have 
$$(\xi|\xi)_A(y)=\sum_{x\in U_0\cap \cG_y}e(x,y),$$
which is continuous and $(\xi|\xi)_A(y_0)=e(x_0,y_0)\geq 1$. 
Note that by assumption $\cG_y\neq \emptyset$ for each $y\in V_1$, 
and so $U_0\cap\ \cG_y\neq \emptyset$ for $y$ sufficiently close to $y_0$. 
Replacing $U_0$ and $V_1$ with smaller neighborhoods 
$U$ of $x_0$ and $V$ of $y_0$ 
respectively, we get the first half of the statement. 
Since $e$ is a positive function, for all $(x,y)\in (U\times V)\cap \cG$ 
we have
$e(x_0,y_0)+\varepsilon\geq e(x,y),$
which means that $e$ is upper semicontinuous. 
\end{proof}

Since $e$ is an upper semicontinuous function on a compact set $\cG$, 
$e$ is bounded above. 
Therefore $(\xi|\eta)_A$ gives an $A$-valued inner product of $X$ with 
the corresponding norm equivalent to the uniform norm of $C(\cG)$, and so 
$X$ is a Hilbert bimodule over $A$. 

Following \cite{KW3}, for $f\in C(K)$ we set
$\tf(y)=\sum_{x\in \cG_y}f(x),$ which is not necessarily continuous. 

Let $u_0,u_1,\cdots, u_n\in X$. 
Then for $\xi\in X$ we have 
\begin{eqnarray*}
\sum_{a=0}^n(\xi|\theta_{u_a,u_a}\xi)_A(y)
&=&\sum_{a=0}^n|(u_a|\xi)(y)|^2\\
&=&\sum_{x,x'\in \cG_y}\sum_{a=0}^n 
e(x',y)e(x,y)u_a(x',y)\overline{u_a(x,y)}\xi(x,y)
\overline{\xi(x',y)}.
\end{eqnarray*}
Let 
$$A(y)_{x',x}=\sum_{a=0}^n \sqrt{e(x',y)e(x,y)}u_a(x',y)
\overline{u_a(x,y)}.$$
Then
\begin{equation}\sum_{a=0}^n\theta_{u_a,u_a}\leq I.
\end{equation}
is equivalent to that the matrix $(A(y)_{x',x})_{x',x}$ 
is dominated by 1 for any $y$. 
In particular, (4.1) implies 
$$A(y)_{x,x}=\sum_{a=0}^n e(x,y)|u_a(x,y)|^2\leq 1,
\quad \forall x\in \cG_y,$$
and so 
\begin{equation}
\sum_{a=0}^n(u_a|f\cdot u_a)(y)
=\sum_{a=0}^n \sum_{x\in \cG_y}f(x)e(x,y)|u_a(x,y)|^2
\leq \tf(y)
\end{equation}
for all positive $f\in C(K)$. 

In what follows, we often identify an element in $C(K)^*_+$ with 
the corresponding measure on $K$. 
In our setting, the Perron-Frobenius type operator $F$ is a positive map 
$F : C(K)^* \rightarrow C(K)^* $ given by   
\begin{equation}
F(\mu)(f)=\sup_{u_a}\sum_{a=0}^n\mu((u_a|f\cdot u_a)),\quad 
\mu\in C(K)^*_+,\; f\in C(K)_+,
\end{equation}
where the supremum is taken over all $u_0,u_1, \cdots,u_n\in X$ 
satisfying (4.1). 
We have $||F||=N$. 

Since $K$ is separable, $X$ is countably generated and has a countable 
basis $\{v_i\}_{i=0}^\infty$. 
We introduce an increasing sequence $\{\varphi_n\}_{n=0}^\infty$ of 
positive maps of $C(K)$ by
$$\varphi_n(f)=\sum_{i=0}^n 
(v_i|f\cdot v_i)_A,\quad f\in C(K).$$
Then $F$ is given by the limit 
\begin{equation}F(\mu)(f)=\lim_{n\to\infty}\mu(\varphi_n(f)). 
\end{equation}

\begin{thm} Let $F$ be as above. Then 
$$F(\mu)(f)=\int_K\tf(y)d\mu(y),\quad f\in C(K)_+,\; \mu\in C(K)^*_+.$$ 
\label{thm;explicitPF}
\end{thm}

\begin{proof} 
First, we show the statement for any Dirac measure 
$\mu=\delta_{y}$ on  $y\in K$.  
Note that the above argument implies 
$$F(\mu)(f)\leq \int_K\tf(y)d\mu(y),\quad f\in C(K)_+,\; \mu\in C(K)^*_+.$$

We fix $y_0\in K$ and $\varepsilon >0$. 
Let $\cG_{y_0}=\{x_0,x_1,\cdots, x_n\}$. 
Applying Lemma \ref{lemma;openneibourhood}, we choose open neighbourhoods 
$U_a$ of $x_a$ and $V$ of $y_0 $ such that 
$U_a\cap U_b=\emptyset$ for $a\neq b$ and 
$$\bigcup_{a=1}^nU_a\supset \cG_y,\quad \forall y\in V,$$
$$|e(x_a,y_0)-\sum_{x\in U_a\cap V}e(x,y)|<\varepsilon , 
\quad \forall y\in V.$$ 
Let $u_a\in X$ be a continuous non-negative function such that 
$\mathrm{supp}(u_a)\subset U_a\times V$ and 
$$u_a(x_a,y_0)=\frac{1}{\sqrt{e(x_a,y_0)+\varepsilon}},$$
$$0\leq u_a(x,y)\leq \frac{1}{\sqrt{e(x_a,y_0)+\varepsilon}},\quad \forall 
(x.y)\in (U_a\times V)\cap \cG.$$
For $y\in V$, we define the matrix $A(y)=(A(y)_{x',x})_{x',x\in\cG_y}$ 
as before. 
We set  $I_a(y)=U_a\cap \cG_y$ for $y\in V$ and 
$A_a(y)=(A(y)_{x',x})_{x',x\in I(a)}$.
Then $A(y)$ is a direct sum of the $A_a(y)$ s. 
Let $c\in \ell^2(I(a))$. 
Then 
\begin{eqnarray*}
\inpr{A_a(y)c}{c}&=&
\sum_{x',x\in I_a(y)}\sum_{b=0}^n 
\sqrt{e(x',y)e(x,y)}u_b(x',y)\overline{c(x')}\overline{u_b(x,y)}c(x)\\
&=&|\sum_{x\in I_a(y)}\sqrt{e(x,y)} \;
\overline{u_a(x,y)}c(x)|^2\\
&\leq &||c||^2\sum_{x\in I_a(y)}e(x,y)|u_a(x,y)|^2\\
&\leq &||c||^2\sum_{x\in I_a(y)}\frac{e(x,y)}{e(x_0,y_0)+\varepsilon} 
 \leq ||c||^2,
\end{eqnarray*}
where we use Lemma \ref{lemma;openneibourhood}. 
Thus $A_a(y)\leq 1$ and (4.1) is satisfied. 
Since 
\begin{eqnarray*}
\sum_{a=0}^n(u_a|f\cdot u_a)_A(y_0)
&=&\sum_{x\in \cG_{y_0}}\sum_{a=0}^n f(x)e(x,y_0)|u_a(x,y_0)|^2 \\
&=&\sum_{x\in \cG_{y_0}}f(x)\frac{e(x,y_0)}{e(x,y_0)+\varepsilon}
\geq \tf(y)\frac{1}{1+\varepsilon},
\end{eqnarray*}
we get $F(\delta_{y_0})(f)=\tf(y_0)$. 

The above computation implies that for every $f\in C(K)_+$ and every $y\in K$, 
$\{\phi_n(f)(y)\}_{n=0}^\infty$ increasingly converges to $\tf(y)$. 
Thus thanks to the bounded (or monotone) convergence theorem, we get 
$$F(\mu)(f)=\lim_{n\to\infty}\int_K\phi_n(f)(y)d\mu(y)
=\int_K\tf(y)d\mu(y),\quad \forall \mu\in C(K)^*_+,$$
which finishes the proof in the general case. 
\end{proof}
 
We denote by $D(e)$  the set of discontinuous points of $e$. 
For example, if $R$ is a rational function, $\cG = \graph R$ and 
$e(z,R(z)) = e(z)$ is the branch index, 
then $D(e)$ is the set $\{(z,R(z)) \in \graph R |\; z \in \cB(R) \}$ 
and $p_1(D(e)) = \cB(R)$, where  $\cB(R)$ is the set of 
branched points of $R$. 
If $\gamma = (\gamma_1,\dots , \gamma_N)$ is a system of    
proper contractions on $K$,  
${\mathcal G} =  \cup _{i=1}^N \{(x,y) \in K^2 ; x = \gamma _i(y)\}$, 
and $e(x,y) =\# \{i | x = \gamma _i(y)\}$, 
then $D(e) = \{(x,y) \in K^2 | x = \gamma _j(y) = \gamma _{j'}(y) 
\text{ for some } j \ne j'\}$ 
and $p_1(D(e)) = \cB(\gamma)$. 
 
\begin{lemma} Let $D(e)$ be the set of discontinuous points of $e$. 
\begin{itemize}
\item [$(1)$] For every $(x_0,y_0)\in \cG\setminus D(e)$, 
there exist closed neighborhoods  $F$ of $x_0$ and $G$ of $y_0$, and a 
continuous map $\gamma:F\rightarrow G$ such that 
$$(F\times G)\cap \cG=\{(\gamma(y),y)\in K^2;\; y\in G\}.$$ 
\item [(2)] $D(e)$ is a closed set. 
\end{itemize}
\label{lemma;discontinuous}
\end{lemma}

\begin{proof} (1) 
Let $0<\varepsilon<1/2$ and take a neighborhood $W$ of 
$(x_0,y_0)$ such that 
$$|e(x_0,y_0)-e(x,y)|<\varepsilon,\quad \forall (x,y)\in W.$$
By Lemma \ref{lemma;openneibourhood}, 
there exist closed neighborhoods $F$ of $x_0$ and $G$ of 
$y_0$ such that $F\times G\subset W$, 
$p_2((F\times G)\cap \cG)=G$, and 
$$|e(x_0,y_0)-\sum_{x\in \cG_y\cap F}e(x,y)|<\varepsilon,
\quad \forall y\in G.$$ 
Since $e(x,y)\geq 1$ for all $(x,y)\in \cG$, the set 
$\cG_y\cap F$ is a singleton for all $y\in G$. 
Thus there exists a function $\gamma:G\rightarrow F$ such that 
$$(F\times G)\cap \cG=\{(\gamma(y),y);\; y\in G\}.$$ 
Since $(F\times G)\cap \cG$ is closed, $\gamma$ is continuous. 

(2) We choose $\xi\in X$ such that $\mathrm{supp}(\xi)\subset F\times G$, 
$0\leq \xi\leq 1$, and $\xi(x,y)=1$ in a neighborhood $W_1$ of 
$(x_0,y_0)$. 
Then for $(x,y)\in W_1$, we have 
$$(\xi|\xi)_A(y)=\sum_{y\in F\cap \cG_y}e(x,y)|\xi(x,y)|^2=e(\gamma(y),y),$$
which shows that $e$ is continuous on $W_1\cap \cG$. 
Thus $D(e)$ is a closed set.  
\end{proof}

The following Proposition shows that the set of 
singularity corresponds to the ideal $I_X$ of $A$. 
We can prove it using the preceding lemmas as in 
\cite{KW1} and in \cite{KW2}.  
It is also shown in a general situation by Muhly and Tomforde 
\cite[Corollary 3.12]{MT}. 

\begin{prop} In the above situation 
$$
I_X = \{f \in C(K) ; f|_{p_1(D(e))} = 0 \} 
\cong  C_0(K\setminus p_1(D(e))).
$$  
\label{proposition;I_X}
\end{prop}

\begin{cor} Let the notation be as above. Then  
\begin{itemize}
\item [$(1)$] If $\beta>\log N$, then there is a one-to-one 
correspondence between $\mathrm{ex}(K_\beta(\alpha))_f
=\mathrm{ex}(K_\beta(\alpha))$ and $p_1(D(e))$. 
\item [$(2)$] If $\beta=\log N$, then there is a one-to-one 
correspondence between $K_\beta(\alpha)_i$ and the set of 
Borel probability measures $\mu$ on $K$ satisfying
$$\int_K f(x)d\mu(x)=\frac{1}{N}\int_K\tf(x)d\mu(x),\quad \forall 
f\in C(K).$$ 
In particular, such a measure $\mu$ satisfies 
$\mu\{y\in K;\#\cG_y<N \}=0.$
\end{itemize}
\label{cor;characterization}
\end{cor} 

\begin{proof} (1) follows from Proposition~\ref{proposition;I_X} 
and Lemma~\ref{lemma;finite}. 
The first statement of (2) follows from Theorem~\ref{thm;K1K2}. 
The case with $f=1$ implies the second statement.   
\end{proof}
\section{classification of KMS states in 
the case of rational functions}

Throughout this section we assume that $R$ is a rational function of degree at least two. 
We shall completely classify the KMS states for the gauge actions  on the $C^*$-algebras 
${\mathcal O}_R(\hat{\mathbb C})$ and ${\mathcal O}_R(J_R)$.  
We can recover the degree of $R$, the number of branched points, 
the number of exceptional points and the orbits of exceptional points 
from the information of the structure  of the KMS states. 

Recall that for $f \in C(\hat{\mathbb C})$ and $ y \in \hat{\mathbb C}$, 
we have 
$$\tf(y)=\sum_{x\in R^{-1}(y)}f(x),$$ 
which is not necessarily continuous. 
For a bounded regular Borel measure $\mu$ on $\hat{\mathbb C}$, 
we denote by $\tau _{\mu} \in C(\hat{\mathbb C})^*$ the 
corresponding positive linear functional (though  
we often identify $\mu$ with $\tau _{\mu}$). 
Then the Perron-Frobenius type
operator $F : C(\hat{\mathbb C})^* \rightarrow 
C(\hat{\mathbb C})^* $ with norm $N = \deg R$ is given by 
\[
F(\tau _{\mu})(f) = \int_{\hat{\mathbb C}} \tf(z) d\mu(z).
\]
In particular, for a Dirac measure $\delta _y$ on $y$ we get  
\[
F(\delta _y) = \sum _{x \in R^{-1}(y)}\delta _x .
\] 

In our particular situation, the conditions (K1) and (K2) in 
Theorem \ref{thm;K1K2} take the following forms: 
$$
\frac{1}{e^{\beta}} 
\int_{\hat{\mathbb C}} \ta(z) d\mu(z) 
=\int_{\hat{\mathbb C}} a(z) d\mu(z),\quad  \forall a \in C_0(\hat{\C}\setminus B(R)),  \leqno (K1)   
$$
$$\frac{1}{e^{\beta}} 
\int_{\hat{\mathbb C}} \ta(z) d\mu(z) 
\leq \int_{\hat{\mathbb C}} a(z) d\mu(z),\quad \forall a \in A_+,\leqno (K2) 
$$
where we identify $C_0(\hat{\C}\setminus \cB(R))$ with 
$I_X=\{f \in A | f(z) = 0 \mbox{ for } z \in \cB(R) \}$.    

\begin{lemma}
If a measure $\mu$  satisfies  the 
condition (K1) and (K2), then the point mass 
$\mu(\{w\}) $  for  $w \in \hat{\mathbb C}$ 
must satisfy the following:
$$e^{-\beta}\mu(\{R(w)\}) = \mu(\{w\}),\quad  
w \not\in  \cB(R), \leqno (K1)'$$
$$e^{-\beta}\mu(\{R(w)\}) \leq \mu(\{w\}),\quad w \in \hat{\mathbb C}. 
\leqno (K2)'$$
\label{lemma;pointmassinequality}
\end{lemma}

\begin{proof} Let $w \not\in  \cB(R)$. 
Then we can approximate the characteristic function $\chi_{\{w\}}$ 
by a monotone decreasing sequence $(a_n)_n$ of non-negative 
functions in $I_X$. Thus
\begin{align*}
\mu(\{w\}) & 
= \lim _n \int_{\hat{\mathbb C}} a_n(z) d\mu(z) 
= \lim _n  e^{-\beta} \int_{\hat{\mathbb C}} 
\sum_{x\in R^{-1}(z)}a_n(x) d\mu(z) \\
& = e^{-\beta}  \int_{\hat{\mathbb C}} 
 \sum_{x\in R^{-1}(z)} \chi_{\{w\}}(x) d\mu(z) 
= e^{-\beta}\mu(\{R(w)\}).
\end{align*}
by the Lebesgue convergence theorem. 
The other  one is similarly obtained. 
\end{proof}

Let us recall some facts on exceptional points for $R$. 
For $z$ and $w$ in $\hat{\mathbb C}$, we define $z \sim _R w$ if 
there exists non-negative integers $n$ and $m$ with 
$R^n(z) = R^m(w)$.  Then $\sim _R$ is an equivalence relation on 
$\hat{\mathbb C}$.  We denote the equivalence class containing 
$z$ by $[z]_R$.  We also define the {\it backward orbit} 
$O^-(z)$ of $z$
by 
\[
O^-(z) := \{ w \in \hat{\mathbb C} | R^n(w) = z \text{ for \ some 
\ non-negative }  n \in \mathbb Z
 \} 
\]
Then $O^-(z) \subset [z]_R$.  

\begin{defi} A point $z$ in $\hat{\mathbb C}$ is 
an {\it exceptional point} for $R$ if the backward orbit 
$O^-(z)$ of $z$ is finite. 
We denote by $E_R$ the set of exceptional points. 
\end{defi}

It is known that $E_R$ is a subset of $F_R \cap \cB(R)$ and that $z$ is 
an exceptional point if and only if $[z]_R$ is finite. 
A rational function $R$ of degree at least two has at most 
two exceptional points by the Riemann-Hurwitz formula. 
The reader is referred to \cite[section 4.1]{B} for basic properties of 
exceptional points, including the following:  

\begin{lemma} Let $R$ be a rational function of degree at least two. 
Then either of the following holds: 
\begin{itemize}
\item[$(1)$] $E_R = \emptyset$. 

\item[$(2)$] $E_R$ consists of one point $z_0$ and $\{z_0\} = [z_0]_R$. 
After a suitable conjugation by a M\"{o}bius transformation, we may assume 
that $z_0 = \infty$ and $R$ is a polynomial.

\item[$(3)$] $E_R$ consists of two points $z_0 \not= z_1$ and 
$[z_0]_R = \{z_0\}$ and $[z_1]_R = \{z_1\}$. 
After a suitable conjugation by a M\"{o}bius transformation, we may assume 
$z_0 = 0$, $z_1 = \infty$, and  $R(z) =z^d$ for some 
positive integer $d$. 

\item[$(4)$] $E_R$ consists of two points $z_0 \not= z_1$ and 
$[z_0]_R = [z_1]_R = \{z_0,z_1\}$. 
After a suitable conjugation by a M\"{o}bius transformation, we may assume 
$z_0 = 0$, $z_1 = \infty$, and $R(z) =z^d$ for some negative integer $d$. 
\end{itemize}
\label{lemma;exceptional}
\end{lemma}

The following proposition is crucial for the complete classification 
of the KMS states. 

\begin{prop} 
Assume that a measure $\mu$ satisfies  the condition (K1) and (K2).  
If $\mu$ has a point mass at $z$ for some $z \not\in E_R$, 
then $\beta > \log N$. 
\label{prop;pointmass}  
\end{prop}

\begin{proof} Thanks to Lemma~\ref{lemma;pointmassinequality}, we have 
$\mu\{w\}>0$ for any $w\in O^-(z)$. 
To prove the statement, it suffices to show that there exists 
$w\in O^-(z)$ satisfying the following two conditions: 
(1) $O^-(w)\cap \cC(R)=\emptyset$, (2) $R^{-m}(w)\cap R^{-n}(w)=\emptyset$ 
for all distinct non-negative integers $m,n$. 
Indeed, Lemma~\ref{lemma;pointmassinequality} with such $w$ implies , 
\[
\sum _{n = 0} ^{\infty}
 (\frac{N}{e^{\beta}})^n {\mu}(\{w\}) \leq \mu(O^-(w))
\leq 1. 
\] 

Assume that (1) does not holds for any $w\in O^-(z)$. 
Then since $\cC(R)$ is a finite set, there exists two positive integers 
$m<n$ such that $R^{-m}(z)\cap R^{-n}(z)$ is not empty, and so 
$z=R^{n-m}(z)$. 
Let $p$ be the minimal positive integer such that $R^p(z)=z$ and 
set $L=\{R^k(z)\}_{k=0}^{p-1}$. 
For any $w\in O^-(z)$, the same argument as above shows that there 
exists a positive integer $q$ such that $w=R^q(w)$. 
Thus for a sufficiently large integer $k$, we get 
$w=R^{qk}(w)\in L$. 
Since $O^{-}(z)$ is an infinite set, this is a contradiction. 
Therefore (1) holds for some $z_1\in O^-(z)$. 

Assuming that (2) does not hold for any $w\in O^-(z_1)$, we get 
a contradiction in the same way, and so there exists $w\in O^-(z_1)$ 
satisfying (1) and (2). 
\end{proof}

We recall basic properties of the Lyubich measure now. 

\begin{defi}[Lyubich measure \cite{FLM}, \cite{L}]  
Let $N = \deg R$ and $\delta _x$ be the Dirac measure on $x$ 
for $x \in \hat{\mathbb C}$.  
For any $y \in \hat{\mathbb C} \setminus E_R$ 
and each $n \in {\mathbb N}$, we define 
a probability measure $\mu _n^y$ on the Riemann sphere 
$\hat{\mathbb C}$ by 
$$
\mu _n^y = \sum _{x\in R^{-n}(y)} 
           N^{-n}e(x)e(R(x))\dots e(R^{n-1}(x))  \delta _x \ . 
$$ 
The sequence $(\mu _n^y)_n$ converges 
weakly to a measure $\mu ^L$, which is called 
the {\it Lyubich measure}.  
The measure $\mu ^L$ is independent 
of the choice of $y \in \hat{\mathbb C} \setminus E_R$. 
\end{defi}

The measure $\mu ^L$ has been studied by several authors, e.g. 
Freire-Lopes-M\"{a}n\'{e} \cite{FLM} and Lyubich \cite{L}.  
The support of $\mu^L$ is the Julia set
$J_R$ and $\mu^L$ is an invariant measure in the sense that 
$\mu ^L(E) = \mu ^L(R^{-1}(E))$ for any Borel set $E$. 
Hence $\int f(R(x)) d\mu ^L(x) = \int f(x) d\mu ^L(x)$. Moreover 
$\mu ^L(R(E)) = N\mu(E)$ for any Borel set $E$ on which 
$R$ is injective. The measure $\mu^L$ is the unique measure of 
maximal entropy. 

We introduce an operator 
$G: C(\hat{\mathbb C}) \rightarrow C(\hat{\mathbb C})$ by setting   
\[
 G(f)(w) = \sum_{z \in R^{-1}(w)} e(z)f(z)
\] 
and the contraction $\overline{G} = N^{-1}G$. 
The conjugate opeator $G^*$ is 
slightly different from $F$. But, 
if $\mu$ has no point mass on  $\cC(R)$, then 
$G^*$ satisfies  $G^* \tau_{\mu} = F\tau_{\mu}$.

\begin{prop}[Lyubich \cite{L}]. 
For any $a \in C(\hat{\mathbb C})$ and compact subset 
$K \subset \hat{\mathbb C}$ with $K \cap E_R = \emptyset$, 
we have  
\[
  \lim_{m \to \infty} \|\overline{G}^m a - \int_{\hat{\mathbb C}} a(z) d\mu^L(z)
  \|_K  =  0, 
 \]
where $\Vert \ \Vert_K$ is the uniform norm on $K$.  
\label{prop:Lyu}
\end{prop}

We use the symbol $\tau ^L$ for $\tau_{\mu^L}$  for simplicity. 

\begin{prop}
Let $\mu$ be a bounded regular Borel measure on $\hat{\mathbb C}$. 
Suppose that $\mu$ satisfies one of the following conditions: 
\begin{itemize}
\item[(1)] $\mu$ satisfies the condition  (K1) and (K2) and 
$\mu$ has no point mass on $\cB(R)$. 
\item [(2)] $\mu$ satisfies the condition  (K1) and 
$\mu$ has no point mass on $\cB(R) \cup \cC(R)$.
\end{itemize}
Then we have the following:  If $\beta \not = \log N$, then $\mu = 0$.
If $\beta = \log N$, then $\mu$ coincides with 
the Lyubich measure up to constant. 
\label{prop;Lyubich}
\end{prop}

\begin{proof} Suppose that $\mu$ is non-zero.  
We may and do assume that 
$\mu$ is a probability measure. 
First we consider the case (1). Since  
$e^{-\beta}\mu(\{R(w)\}) \leq \mu(\{w\})$ for $w \in \hat{\mathbb C}$ 
by Lemma \ref{lemma;pointmassinequality}, 
$\mu$ has no point mass on $\cC(R)$. 
Therefore (K1) holds for all $a \in A$ and 
we may consider the case (2). 
For $a = 1$ , we have $\ta(x) = N $ a.e. $x$ with respect to 
$\mu$ and $\beta = \log N$.  
For any $a \in A$, we have 
$$\tau_{\mu}(a) = F_{\beta}(\tau_{\mu})(a) = 
\tau_{\mu}(\overline{G}(a)).$$ 
Note that $\mu(E_R) = 0$ holds as $E_R \subset \cB(R)$. 
For any $\varepsilon > 0$, we can choose a  compact subset 
$K \subset \hat{\mathbb C}$ satisfying $K \cap E_R = \emptyset$ 
and $\mu(K^c) \|a\| < \varepsilon $ by the regularity of $\mu$. 
Thus   
\begin{eqnarray*}\lefteqn{
  | \tau_{\mu}(a) - \tau^L(a)| 
 = |\tau_{\mu}(\overline{G}^m(a)) - \tau^L(a)|} \\
 &= & | \int_K \overline{G}^m(a)\,d\mu
    + \int_{K^c} \overline{G}^m(a)\,d\mu
    - (\tau^L(a) \mu(K) + \tau^L(a) \mu(K^c))| \\
 &\le &  | \int_K \big(\overline{G}^m(a) - \tau^L(a)\big) \,d\mu|
     + | \int_{K^c} \overline{G}^m(a)\,d\mu | 
+  \tau^L(a) \mu(K^c)
 < 3 \varepsilon. 
\end{eqnarray*} for sufficiently large $m$ thanks to 
Proposition \ref{prop:Lyu}, which shows $\mu = \mu^L$.  
\end{proof}

\begin{prop} The Lyubich measure $\mu^L$ 
is extended to an infinite type $\log N$-KMS state 
$\varphi^L$ on ${\mathcal O}_R(\hat{\mathbb C})$ for 
the gauge action. 
\end{prop}

\begin{proof}
Since the  Lyubich measure has no point mass, 
it satisfies $F(\tau^L)=G^*(\tau^L)$. 
On the other hand, for any $x\in \hat{\C}$ we have 
$$G^*(\mu_n^x)=N\sum_{y\in R^{-1}(x)}\mu_{n+1}^y,$$
which implies $G^*(\tau^L)=N\tau^L$.   
Thus by Corollary \ref{cor;characterization}, 
the Lyubich measure $\mu ^L$ corresponds to a 
$\log N$-KMS state $\varphi^L$.
\end{proof}

Now we are ready to state our classification results. 
Since we already know the result for $\beta>\log N$ thanks to Corollary 
\ref{cor;characterization}, we assume $0<\beta\le \log N$.  
Proposition~\ref{prop;pointmass} implies that every finite type 
$\beta$-KMS state, if it exists, arises from a point in $E_R$. 
Let $\mu$ be a probability measure corresponding to an extreme 
infinite type $\beta$-KMS state $\varphi$, which satisfies $F_\beta(\mu)=\mu$. 
Let $\mu=\mu_a+\mu_d$ be the decomposition of $\mu$ into the 
atomic part $\mu_a$ and the diffuse part $\mu_d$. 
Then a similar argument as in the proof of 
Lemma~\ref{lemma;pointmassinequality} shows that 
$F_\beta(\mu_a)$ is atomic and $F_\beta(\mu_d)$ is diffuse again, 
and so we get $F_\beta(\mu_a)=\mu_a$ and $F_\beta(\mu_d)=\mu_d$. 
Therefore either $\mu_a=0$ or $\mu_d=0$ holds. 
If $\mu_d=0$, Proposition~\ref{prop;pointmass} and 
Lemma~\ref{lemma;exceptional} imply that $\mu$ would be supported by  
$E_R$. 
Since $R^{-1}(x)$ is a singleton for $x\in E_R$, this does not occur. 
Therefore $\mu_a=0$ and Proposition \ref{prop;Lyubich} implies that 
$\beta=\log N$ and $\mu$ is the Lyubich measure. 

The above observation with an easy case-by-case analysis using 
Lemma~\ref{lemma;exceptional} shows the following theorems:

\begin{thm} 
For $\beta>0$, there exists an infinite type $\beta$-KMS state for 
the gauge action $\alpha$ on $\cO_R(\hat{\C})$ if and only if 
$\beta=\log N$. 
The set $K_{\log N}(\alpha)_i$ is a singleton consisting of 
$\varphi_L$ given by the Lyubich measure. 
\end{thm} 

\begin{thm} For the finite type $\beta$-KMS states for the gauge action 
$\alpha$ on $\cO_R(\hat{\C})$, the following holds:
\begin{itemize}
\item [$(1)$] 
For $\beta>\log N$, the set  
$\mathrm{ex}(K_\beta(\alpha)_f)$ of the extreme finite type $\beta$-KMS states 
$\varphi$ is parameterized by the branched points $\cB(R)$. 
For $w \in \cB(R)$, the corresponding extreme $\beta$-KMS state 
$\varphi _{\beta,w}$ is determined by 
its restriction $\varphi _{\beta,w} |_{C(\hat{\mathbb C})}$, 
i.e.,  the regular Borel measure $\mu_{\beta,w}$ on $\hat{\mathbb C}$ 
given by 
\[
\mu_{\beta,w} = m_{\mu _{\beta,w}}\sum_{k=0}^{\infty} \frac{1}{e^{k\beta}} 
 \sum_{z \in R^{-k}(w)} \delta_z, 
\]
where $m_{\mu _{\beta,w}}$ is the normalized constant.  

\item [$(2)$] For $0<\beta\leq \log N$, the set 
$\mathrm{ex}(K_\beta(\alpha)_f)$ of the extreme finite type $\beta$-KMS states 
coincides with $\{\varphi_{\beta,w}\}_{w\in E_R}$.  
\end{itemize}
\label{thm;classification}
\end{thm} 

Furthermore the set of all extreme $\beta$-KMS states, which depends only  
on the system $(\cO_R(\hat{\C}),\alpha)$,  is obtained just as 
$\mathrm{ex}(K_\beta(\alpha))=\mathrm{ex}(K_\beta(\alpha)_f)\cup 
\mathrm{ex}(K_\beta(\alpha)_i)$.

\begin{rema} 
Corollary \ref{cor;characterization} shows that the GNS representation of 
every extreme finite type KMS state with $\beta>0$ is a type I factor representation. 
In Section 7, we show that the GNS representation of $\varphi^L$ is a factor 
representation of type III$_{1/N}$.   
\end{rema}

\begin{rema}
Theorem~\ref{thm;K1K2} still holds for $\beta=0$ if we  
define a 0-KMS state to be an $\alpha$-invariant tracial state.  
Such a state exists if and only if $E_R$ is not empty and we have 
the following: 
\begin{itemize}
\item [(1)] When $E_R$ consists of one point $w$, then there exists a unique  
$\alpha$-invariant trace state $\varphi_w$. 
The restriction of $\varphi_w$ to $C(\hat{\C})$ is given by 
the Dirac measure $\delta_w$.  
\item [(2)] When $E_R$ consists of two points $w_1$ and $w_2$ with $R(w_1)=w_1$ 
and $R(w_2)=w_2$, the set of $\alpha$-invariant trace states has exactly two 
extreme points $\{\varphi_{w_i}\}_{i=1,2}$. 
The restriction of $\varphi_{w_i}$ to $C(\hat{\C})$ is given by 
$\delta_{w_i}$ for $i=1,2$. 
\item [(3)]  When $E_R$ consists of two points $w_1$ and $w_2$ with $R(w_1)=w_2$ 
and  $R(w_2)=w_2$, then there exists a unique $\alpha$-invariant
trace $\varphi$. 
The restriction of $\varphi$ to $C(\hat{\C})$ is given by 
$$\frac{1}{2}(\delta_{w_1}+\delta_{w_2})$$
\end{itemize} 
Note that the GNS representations of these states are not factor representations. 
It is routine work to show that they give finite type I von Neumann algebras. 
\end{rema}

\begin{exam} Let $R(z)=z^N$ with $N\geq 2$. 
Then $E_R=\{0,\infty\}=\cB(R)$ with $R(0)=0$ and $R(\infty)=\infty$. 
For every $\beta>0$ and $w=0,\infty$, we have $\mu_{\beta,w}=\delta_w$.   
\end{exam} 

\begin{exam} Let $R(z)=z^{-N}$ with $N\ge 2$. 
Then we have $E_R=\cB(R)=\{0,\infty\}$ with $R(0)=\infty$ and $R(\infty)=0$. 
For every $\beta>0$ we have 
\[
 \mu_{\beta,0} = \frac{e^\beta}{e^{\beta} + 1} \delta _{0}
 + \frac{1}{e^{\beta} + 1} \delta _{\infty},
\]
\[
 \mu_{\beta,\infty} = \frac{1}{e^{\beta}+ 1} \delta _{0}
 + \frac{e^\beta}{e^{\beta} + 1} \delta _{\infty}. 
\] 
\end{exam}

\begin{exam} Let $R(z) = z^N +1$ with $N \ge 2$. Then 
$E_R = \{\infty\}$ and $\cB(R) = \{0,\infty\}$. 
For any $\beta \geq 0$ the Dirac measure $\delta_{\infty}$ extends to a 
$\beta$-KMS state $\varphi_{\beta,\infty}$. 
The measure  
\[
\mu_{\beta,0} =(1-\frac{N}{e^\beta})\sum_{k=0}^{\infty} (\frac{1}{e^{k\beta}} 
 \sum_{z \in R^{-k}(0)} \delta_z), 
\]
corresponds to a $\beta$-KMS state $\varphi_{\beta,0}$ for $\beta > \log N$. 
For $0 \le \beta < \log N$, the set of $\beta$-KMS states consists of one 
point $\varphi_{\beta,\infty}$. 
For $\beta =\log N$, the set of extreme  $\beta$-KMS states 
consists of two points $\varphi^L$ and $\varphi_{\beta,\infty}$. 
For $\log N < \beta $, the set of extreme  $\beta$-KMS states 
consists of two points $\varphi_{\beta,\infty}$ and $\varphi_{\beta,0}$. 
\end{exam}

Finally we consider the $C^*$-algebras ${\mathcal O}_R = {\mathcal O}_R(J_R)$ 
associated with a rational function $R$ on the Julia set $J_R$, 
which is purely infinite and simple \cite{KW1}. 

\begin{thm} 
For the gauge action $\alpha$ on the $C^*$-algebras 
${\mathcal O}_R(J_R)$, the following hold:
\begin{itemize}
\item [$(1)$] For $0 < \beta < \log N$, there exist no $\beta$-KMS states. 
\item [$(2)$] For $\beta = \log N$, there exists a unique $\beta$-KMS state 
$\varphi^L$ and its restriction $\varphi^L |_{C(J_R)}$ is the Lyubich measure. 
The GNS representation of $\varphi^L$ is a factor representation of 
type $III_{1/N}$.  
\item [$(3)$] For $\log N < \beta$, the set $\mathrm{ex}(K_\beta(\alpha))$ 
of extreme $\beta$-KMS states is parameterized by the branched points 
$\cB(R) \cap J_R$. 
For $w \in \cB(R) \cap J_R$, the corresponding extreme $\beta$-KMS state 
$\varphi _{\beta,w}$ is determined by its restriction 
$\varphi _{\beta,w} |_{C(J_R)}$, which is given 
by the same formula as in Theorem \ref{thm;classification}. 
The GNS representation of $\varphi_{\beta,w}$ is a factor representation of 
type $I$.
\end{itemize}
\end{thm}

\begin{proof} The proof follows from a similar argument as in the 
case of $\cO_R(\hat{\C})$ with the fact that the Julia set $J_R$ contains 
no exceptional point by \cite{B}.  
\end{proof}

\section{classification of KMS states in the case of self-similar sets}
Let $(K,d)$ be a compact metric space and let 
$\gamma=(\gamma_1,\gamma_2,\cdots,\gamma_N)$ be a system of proper 
contractions satisfying $K=\bigcup_{i=1}^N\gamma_i(K)$. 
In this section, we consider the structure of the KMS states for the 
gauge action on $\cO_X=\cO_\gamma (K)$. 
Recall $A = C(K)$ and $I_X = \{f \in C(K) ; f|_{\cB(\gamma)} = 0 \} 
=C_0(K\setminus \cB(\gamma))$ in this case. 
The set  $\cB(\gamma)$ is finite if and only if $\cC(\gamma)$ is finite. 

For $y\in K$, we set $\gamma(y)=\bigcup_{i=1}^N\{\gamma_i(y)\}$.  
The Borel function $\tilde{a}$ for $a \in A$ is 
\begin{equation}
 \tilde{a}(y) = \sum_{x \in \gamma(y)}a(x)
 = \sum_{j=1}^N \frac{1}{e(\gamma_j(y),y)}a(\gamma_j(y)).  
\end{equation}
Note that if $\cC(\gamma)$ is not empty, $\tilde{a}$ is not necessarily 
continuous.  
Theorem \ref{thm;explicitPF} shows  
\begin{equation}
F_\beta(\delta_y)=e^{-\beta}\sum_{j=1}^N
\frac{1}{e(\gamma_j(y),y)}\delta_{\gamma_j(y)} 
= e^{-\beta}\sum_{x \in \gamma(y)}\delta_{x}.
\end{equation}
If a probability measure $\mu$ on $K$ satisfies (K2), we have 
$\mu\geq F_\beta(\mu)\geq \mu(\{x\})F_\beta(\delta_x)$ and so, 
the following holds: 
\begin{equation}
\mu(\{\gamma_i(x)\})\geq e^{-\beta}\mu(\{x\}).
\label{pointmass}
\end{equation}



\begin{lemma} Let $\mu$ be a probability regular Borel  measure 
on $K$ satisfying (K1) and (K2). If $\mu(\cC(\gamma)) = 0$, 
then $\log N \leq \beta$. 
\label{lemma;pointmassless}
\end{lemma}
\begin{proof}
Suppose $\mu(\cC(\gamma))=0$. Then 
\[
\tilde{1}(y) 
= \sum_{j=1}^N \frac{1}{e(\gamma_j(y),y)}1(\gamma_j(y)) 
= N, \quad  \text{a.e. $y$ with respect to  $\mu$}.
\] 
(K2) implies that 
\[
\frac{N}{e^{\beta}} =  
\frac{1}{e^{\beta}} 
\int_{K} \tilde{1}(z) d\mu(z) 
\leq \int_{K} 1 d\mu(z)  
= 1 .
\]
Hence $\log N \leq \beta$.
\end{proof}

For $x\in K$ and $n\geq 0$, we set  the $n$-th orbit of $x$ by  
$$O_n(x)=\{\gamma_{i_1}\cdots\gamma_{i_2}\cdots\gamma_{i_n}(x)\in K;\: 
1\leq i_1,i_2,\cdots,i_n\leq N\}$$
and set $O(x)=\bigcup_{n=0}^\infty O_n(x)$. 
Note that if a measure $\mu$ satisfying (K2) has a point mass at $x$, 
it has a point mass at $y$ for all $y\in O(x)$ thanks to (\ref{pointmass}). 

\begin{lemma} Let $0<\beta\leq\log N$ and 
let $\mu$ be a probability regular Borel measure 
on $K$ satisfying (K1) and (K2). 
Assume either $0<\beta<\log N$ or $\mu$ is of finite type. 
Let $y\in K$. 
If there exists $x\in O(y)$ such that $O(x)\cap \cC(\gamma)=\emptyset$, 
then $\mu$ has no point mass at $y$.  
\label{lemma;selfsimilarpointmass}
\end{lemma}

\begin{proof} 
First we assume $0 < \beta<\log N$. 
If $\mu$ had a point mass at $y$, then 
$\mu$ would have a point mass at $x$ as well and 
we would have
$\mu\geq \mu(\{x\})F_\beta^n(\delta_x)$. 
However, by assumption we get 
\[
1 \geq \mu(\{x\})F_\beta^n(\delta _x)(K)
=\mu(\{x\}) (Ne^{-\beta})^n \rightarrow \infty,
\quad (n\to \infty), 
\]
which is a contradiction. 

Now we assume $\beta=\log N$ and the corresponding KMS state is of finite 
type. 
Then there is a measure $\nu$ supported by $\cB(\gamma)$ such that 
$$\mu=\sum_{n=0}^\infty F_\beta^n(\nu).$$
Suppose that $\mu$ has a point mass at $y$. 
Then there would exist $n$ and $c>0$ such that $F_\beta^n(\nu)\geq c\delta_x$, 
and so $F_\beta^{n+m}(\nu)\geq cF_\beta^m(\delta_x)$. 
As we have $F_\beta^m(\delta_x)(1)=1$, 
$$
1 = \mu(1) =\sum_{k=0}^\infty F_\beta^k(\nu)(1) 
\geq \sum_{k=n}^\infty c = \infty.  
$$ 
We get a contradiction. 
\end{proof}
 
To describe infinite type KMS states, we recall the notion of 
the Hutchinson measure for a self-similar set (see \cite{H} for details). 

\begin{defi} Let $\overline{G}:C(K):\rightarrow G(K)$ be the unital positive 
map defined by 
$$\overline{G}(a)=\frac{1}{N}\sum_{i=1}^Na\cdot\gamma_i,\quad a\in C(K).$$
The {\it Hutchinson measure} is a unique $\overline{G}^*$-invariant probability measure. 
\end{defi}

\begin{lemma} Unless $K$ consists of one point, 
the Hutchinson measure $\mu^H$ has no point mass. 
\label{lemma;nopintmass}
\end{lemma}

\begin{proof} Let $\mu^H=\mu_a+\mu_d$ be the decomposition of $\mu^H$ 
into the atomic part $\mu_a$ and the diffuse part $\mu_d$. 
Then it is easy to show that $\overline{G}^*(\mu_a)$ is atomic and 
$\overline{G}^*(\mu_d)$ is diffuse. 
Thus either $\mu_a=0$ or $\mu_d=0$ holds thanks to the uniqueness of a 
$\overline{G}$-invariant measure. 

Suppose $\mu_d=0$. 
We set $m=\max\{\mu\{x\};\;x\in K\}$ and $L=\{x\in K;\; \mu\{x\}=m\}$, 
which is a finite set. 
Since 
$$\overline{G}^*(\delta_x)=\frac{1}{N}\sum_{i=1}^N\delta_{\gamma_i(x)},$$
and $\overline{G}^*(\mu^H)=\mu^H$, we have 
\begin{eqnarray*} m(\#L)&=&\mu^H(L)=\overline{G}^*(\mu^H)(L)
=\frac{1}{N}\sum_{i=1}^N\sum_{y\in K}\mu\{y\}\delta_{\gamma_i(y)}(L)\\
&=&\frac{1}{N}\sum_{i=1}^N\sum_{x\in L,y\in K}\mu\{y\}\delta_{x,\gamma_i(y)}. 
\end{eqnarray*}
This shows that $\gamma_i(L)=L$ for all $i$ and $L=K$. 
Since $\gamma_i$ is a proper contraction and $L$ is a finite set, 
we conclude that $K$ consists of one point.  
\end{proof}

Note that when a measure $\mu$ satisfies $\mu(\cC(\gamma))=0$, we 
have $F(\mu)=N\overline{G}(\mu)$. 
Thus when the Hutchinson measure $\mu^H$ satisfies $\mu^H(\cC(\gamma))=0$, 
it gives rise to an infinite type $\log N$-KMS state, which we denote by 
$\varphi^H$.  
Thanks to Lemma~\ref{lemma;nopintmass}, the condition $\mu^H(\cC(\gamma))=0$ 
holds whenever $\cC(\gamma)$ is a countable set.  
Conversely, we have 

\begin{lemma} Let $\varphi$ be a $\beta$-KMS state of 
infinite type for $\beta=\log N$ and let $\mu$ be the corresponding 
measure on $K$. 
Then $\varphi=\varphi^H$ and $\mu^H(\cC(\gamma))=0$. 
\label{lemma;Hutchinson}
\end{lemma}

\begin{proof} $F_\beta(\mu)=\mu$ implies 
\begin{eqnarray*}
0&=&\int_K(N-\tilde{1}(x))d\mu(x)=\int_K \sum_{j=1}^N
\Big(1-\frac{1}{e(\gamma_j(x),x)}\Big)d\mu(x)\\
&=&\int_{\cC(\gamma)} \sum_{j=1}^N
\Big(1-\frac{1}{e(\gamma_j(x),x)}\Big)d\mu(x),
\end{eqnarray*}
which shows $\mu(\cC(\gamma))=0$. 
Thus $\mu=F_\beta(\mu)=\overline{G}(\mu)$ and so $\mu=\mu^H$.
\end{proof}

\begin{thm} 
Let $\gamma=(\gamma_1,\gamma_2,\cdots,\gamma_N)$ 
be a system of proper contractions on a compact metric space $K$ 
such that $K=\bigcup_{i=1}^{N} \gamma_{i}(K)$.
Suppose either that $\cB(\gamma)$ is empty or that 
$\cB(\gamma)$ is a finite set and for every 
$y\in \cC(\gamma)$, there exists $x\in O(y)$ such that 
$O(x)\cap \cC(\gamma)=\emptyset$. 
For the gauge action $\alpha$ on the $C^*$-algebras 
${\mathcal O}_{\gamma}(K)$, the following hold:
\begin{itemize}
\item [$(1)$] For $0<\beta < \log N$, there exist no $\beta$-KMS states. 

\item [$(2)$] For $\beta = \log N$, then there exists a unique 
$\beta$-KMS state $\varphi^H$ and its restriction $\varphi^H |_{C(K)}$ is 
the Hutchinson measure $\mu^H$. 

\item [(3)] For $\log N < \beta$, the set $\mathrm{ex}(K_\beta(\alpha))$ 
of extreme $\beta$-KMS states $\varphi$
is parameterized by the branched points $\cB(\gamma)$. 
\end{itemize}
\label{thm;selfsimilarkms}
\end{thm}

\begin{proof} 
For $0< \beta<\log N$, Lemma \ref{lemma;pointmassless} and 
Lemma \ref{lemma;selfsimilarpointmass}  show that there is no 
$\beta$-KMS state. 
Assume that $\mu$ is a measure corresponding to a finite type 
$\beta$-KMS state for $\beta=\log N$. 
Then Lemma \ref{lemma;selfsimilarpointmass} implies $\mu(\cC(\gamma))=0$  
and so $\mu(K)=F_\beta(\mu)(K)$ would hold. 
However, this means that $\mu$ is of infinite type, which is a contradiction. 
The rest follows from Lemma \ref{lemma;Hutchinson} and   
Corollary \ref{cor;characterization}. 
\end{proof}

\begin{rema} 
In Section 7, we show that 
the GNS representation of $\varphi^H$ is a factor representation of type 
$III_{1/N}$. 
\end{rema}

\begin{exam}[Inverse branches of the tent map] 
Let $K=[0,1]$, $\gamma_1(y) = \frac{1}{2}y$,  
and $\gamma_2(y) = 1 - \frac{1}{2}y$. 
Then $\cB(\gamma) =\{ \frac{1}{2}\}$ and  $\cC(\gamma) = \{1\}$. 
Since  $O(1/2)\cap \cC(\gamma)=\emptyset$, the assumption of Theorem 
\ref{thm;selfsimilarkms} is satisfied. 
The Hutchinson measure $\mu^H$ is 
the normalized Lebesgue measure on $[0,1]$.  
Hence for $0< \beta < \log 2$, there exist no $\beta$-KMS states. 
For $\beta = \log 2$, there exists a unique $\beta$-KMS state 
$\varphi^H$ and its restriction $\varphi^H |_{C([0,1])}$ is 
the normalized Lebesgue measure. 
For $\log 2 < \beta$, there exists a unique $\beta$-KMS state 
$\varphi_{\beta,1/2}$ and its restriction $\varphi |_{C([0,1])}$ is 
\[
\mu_{\beta,1/2}=(1-\frac{2}{e^\beta})
\sum_{n=0}^{\infty} \frac{1}{e^{n\beta}} 
\sum_{(j_1,\cdots,j_n) \in \{1,2\}^n}
\delta_{\gamma_{j_1}\cdots \gamma_{j_n}(1/2)}.
\]
\end{exam}

\begin{exam} Let $K=[0,1]$, $\gamma_1(y) = \frac{1}{2}y$,  
and $\gamma_2(y) =\frac{1}{2}(y+1)$. 
Then $\cB(\gamma)=\phi$ and $\cC(\gamma)=\phi$ and 
The Hutchinson measure $\mu^H$ is the normalized Lebesgue measure on $[0,1]$.
The $C^*$-algebra ${\mathcal O}_{\gamma}(K)$ with $\alpha$ 
is isomorphic to the Cuntz algebra ${\mathcal O}_2$ with 
the usual gauge action. 
Thus a $\beta$-KMS state on ${\mathcal O}_{\gamma}(K)$ exists if and 
only if $\lambda = \log 2$ and its restriction to 
$C([0,1])$ is the normalized Lebesgue measure.  
\end{exam}

\begin{exam}[Sierpinski Gasket] 
Let $\Omega$ be a regular triangle in ${\mathbb R}^2$ 
with three vertices  
$c_1=(1/2,\sqrt{\,3}/2)$, $c_2 = (0,0)$ and $c_3=(1,0)$.  
The middle point of $c_1c_2$ is denote by $b_1$, the middle point of 
$c_1c_3$ is denoted by $b_2$ and the middle point of $c_2c_3$ is 
denoted by $b_3$.  
We define proper contractions ${\gamma_i}$ for $i=1,2,3$ by 
\[
 {\gamma}_1(x,y)  = \left( \frac{x}{2} + \frac{1}{4}, \frac{y}{2} 
+ \frac{\sqrt{\,3}}{4}\right), \quad
 {\gamma}_2(x,y)  = \left( \frac{x}{2}, \frac{y}{2}\right), \quad 
 {\gamma}_3(x,y)   = \left( \frac{x}{2} + \frac{1}{2},\frac{y}{2} 
 \right).  
\]
Then the self-similar set $K$ associated with $\gamma$ is 
the Sierpinski gasket. 
It is known that $K$ is homeomorphic to 
the Julia set $J_R$ of the rational function 
$$R(z) = \frac{z^3-\frac{16}{27}}{z}.$$
However, these three contractions can not be identified with the 
inverse branches of $R$ because $\gamma _1(c_2) = \gamma _2(c_1)$. 
Let 
$\tilde{\gamma _1} = \gamma _1$, 
$\tilde{\gamma _2} = \rho_{-\frac{2\pi}{3}} \circ \gamma _2$,  
and $\tilde{\gamma _3} = \rho_{\frac{2\pi}{3}} \circ \gamma _3$, 
where $\rho_{\theta}$ is a rotation by the angle $\theta$. 
The self-similar set for 
$\tilde{\gamma}=(\tilde{\gamma _1}, \tilde{\gamma _2},  \tilde{\gamma _3})$ 
is $K$ as well and they are inverse branches of a map $h: K \rightarrow K$, 
which is conjugate to $R : J_R \rightarrow J_R$.  
Therefore ${\mathcal O}_R(J_R)$ is isomorphic to  
${\mathcal O}_{\tilde{\gamma}}(K)$, 
which is known to be a purely infinite and simple $C^*$-algebra. 
The $C^*$-algebras ${\mathcal O}_{\gamma}(K)$ and 
$\cO_{\tilde{\gamma}}(K)$ have non-isomorphic $K$-groups and hence 
they are not isomorphic (See \cite{KW1} for details). 

Since $\cB(\gamma)=\emptyset$ and $\cC(\gamma)=\emptyset$, there exists 
a  $\beta$-KMS state for the gauge automorphism for 
${\mathcal O}_{\gamma}(K)$ if and 
only if $\lambda = \log 3$. 

For $\tilde{\gamma}$, we have $\cB(\tilde{\gamma})=\{b_1,b_2,b_3\}$, 
$\cC(\tilde{\gamma}) =\{ c_1,c_2,c_3 \}$, and 
the assumption of Theorem \ref{thm;selfsimilarkms} is satisfied.  
Thus for $0<\beta < \log 3$, there exist no $\beta$-KMS states. 
For $\beta = \log 3$, there exists a unique $\beta$-KMS state 
whose restriction to $C(K)$ is given by the Hutchinson measure $\mu^H$. 
For $\log 3 < \beta$, the set  
$\mathrm{ex}(K_\beta(\alpha))$ of extreme $\beta$-KMS states is parameterized 
by the branched points $\cB(\gamma) = \{b_1,b_2,b_3\}$.
For $b_i \in \cB(\gamma)$, the corresponding extreme $\beta$-KMS state 
is given by the measure
 \[
\mu_{\beta, b_i} = (1 -\frac{3}{e^{\beta}})
  \sum_{n=0}^{\infty} \frac{1}{e^{n\beta}} 
\sum_{(j_1,\cdots,j_n) \in \{1,2,3\}^n}
  \delta_{\gamma_{j_1}\cdots \gamma_{j_n}(b_i)}. 
\]
\end{exam}

\section{Type III representations by Lyubich and Hutchinson measures}

In this section, we show that the KMS states arising from the Lyubich 
measures and the Hutchinson measures give 
type III$_{1/N}$ representations.  
Our proof is based on the facts that the Cuntz algebra $\cO_N$ arises from 
the Bernoulli $N$-shift and that the unique KMS state of the gauge action 
of $\cO_N$ gives a type III$_{1/N}$ representation. 

Let $N$ be an integer greater than 1 and $K_N=\{1,2,\cdots,N\}^{\N}$.
We denote by $\nu_N$ the Bernoulli measure with weight 
$(\frac{1}{N},\dots,\frac{1}{N})$, i.e. the infinite 
product measure of the uniform measure of $\{1,2,\cdots,N\}$. 
Let $\sigma$ be the Bernoulli shift 
$$\sigma(x_1,x_2,\cdots)=(x_2,x_3,\cdots).$$
For $j=1,2,\cdots,N$, we set $\rho_j$ to be the inverse branches of 
$\sigma$, that is, 
$$\rho_j(x_1,x_2,\cdots)=(j,x_1,x_2,\cdots).$$
Note that $\nu_N$ is nothing but the Hutchinson measure of 
$\rho=(\rho_1,\rho_2,\cdots,\rho_N)$. 

Let $Y$ be the $C(K_N)-C(K_N)$ bimodule giving $\cO_{\rho}(K_N)$. 
Since 
$$\bigcup_{j=1}^N\{(\rho_j(y),y)\in K_N^2;\;x\in K_N\}
=\{(x,\sigma(x))\in K_N^2;\;x\in K_N\},$$ 
in what follows, we identify $Y$ with 
$C(K_N)$ whose Hilbert $C^*$-bimodule structure is given by  
$$(a\cdot f\cdot b)(x)=a(x)f(x)b(\sigma(x))$$
$$(f|g)_{C(K_N)}(y)=\sum_{j=1}^N \overline{f(\rho_j(y))}{g(\rho_j(y))}.$$

We use the symbol $1_E$ for the characteristic function of a subset 
$E\subset K_N$. 
For $\xi=(\xi_1,\xi_2,\cdots,\xi_n)\in \{1,2,\cdots,N\}^n$, 
we denote by $E_\xi$ the cylinder set
$$E_\xi=\{x\in K_N; \;x_k=\xi_k,\; 1\leq k\leq n\}.$$
It is known that the $C^*$-algebra $\cO_\rho(K_N)$ is isomorphic to 
the Cuntz algebra $\cO_N$. 
Indeed, $\cO_\rho(K_N)$ is generated by $S_{1_{E_1}},S_{1_{E_2}},
\cdots,S_{1_{E_N}}$, which satisfy the Cuntz algebra relation 
(see \cite[Proposition 4.1]{KW2} and \cite[Section 4]{PWY}). 
Note that the isomorphism $\cO_N\ni S_j\mapsto S_{1_{E_j}}\in \cO_\rho(K_N)$ 
intertwines the gauge actions, and in consequence, the KMS states. 
We give a $W^*$-version of this argument first. 

\subsection{The case of Hutchinson measures}
Let $\gamma=(\gamma_1,\gamma_2,\cdots,\gamma_N)$ be 
an $N$-tuple of proper contractions of a compact metric space $K$ 
satisfying the self-similar condition 
$K=\bigcup_{j=1}^N\gamma_j(K).$
We assume that the Hutchinson measure $\mu^H$ satisfies 
$\mu_H(\cC(\gamma))=0$ and denote by $\varphi^H$ the corresponding 
KMS state of the gauge action $\alpha$ on $\cO_\gamma(K)$ as before. 

For $\xi=(\xi_1,\xi_2,\cdots,\xi_n)\in \{1,2,\cdots,N\}^n$, we use 
the notations 
$$\gamma_\xi=\gamma_{\xi_1}\cdot\gamma_{\xi_2}\cdots\gamma_{\xi_n}.$$ 
Since all the maps $\gamma_1,\gamma_2,\cdots,\gamma_N$ are proper contractions, 
there exists an surjective continuous map 
$q : K_N \rightarrow K$ such that  
\[
\{q (x)\} =  \bigcap_{n=0}^\infty\gamma_{x_1}\cdot
\gamma_{x_2}\cdots\gamma_{x_n}(K)
\] 
We have 
$q \circ \rho _i = \gamma _i \circ q$ for $i = 1, \dots , N$. 
Since 
$$ 
\sum_{j=1}^N{\gamma_j}_*q_*\nu_N=\sum_{j=1}^Nq_*{\rho_j}_*\nu_N=Nq_*\nu_N,$$
the uniqueness of the Hutchinson measure $\mu^H$ implies $\mu_H= q_* \nu_N$.

\begin{lemma} For all $f,g\in L^2(K,\mu_N)$, the following holds:
$$(f|g)_{\mu}=\lim_{n\to\infty}\frac{1}{N^n}
\sum_{\xi\in \{1,2,\cdots,N\}^n}(f\cdot\gamma_{\xi}|1)_\mu 
(1|g\cdot\gamma_\xi)_\mu,$$
where $(\cdot|\cdot)_\mu$ is the inner product of $L^2(K,\mu)$.  
\end{lemma} 

\begin{proof} Since $q^*:L^2(K,\mu_H)\ni f\mapsto f\cdot q\in L^2(K_N,\nu_N)$ 
is an isometry satisfying $\rho_j^*q^*=q^*\gamma_j^*$, 
it suffices to show the statement for $(K_N,\rho,\nu_N)$ instead of 
$(K,\gamma,\mu)$. 

Let $\cF_n$ be the $\sigma$-field generated by the first $n$ 
coordinate functions of $K_N$, and $E(f|\cF_n)$ be the conditional 
expectation of $f$ given $\cF_n$. 
Then for $f\in L^1(K_N,\nu_N)$ and $\xi\in \{1,2,\cdots,N\}^n$, 
we have 
$$E(f|\cF_n)(\xi)=\int_{K_N}f\cdot\rho_\xi(x)d\nu_N(x).$$
Thus for $f,g\in L^2(K_N,\nu_N)$, we get 
$$\frac{1}{N^n}
\sum_{\xi\in \{1,2,\cdots,N\}^n}
(f\cdot\rho_\xi|1)_{\nu_N}(1|g\cdot\gamma_\xi)_{\nu_N}
=(E(f|\cF_n)|E(g|\cF_n))_{\nu_N},$$
which tends to $(f|g)_{\nu_N}$ as $n$ goes to infinity. 
\end{proof}

Let $A=C(K)$ and $X$ be the $A-A$ bimodule giving $\cO_\gamma(K)$.  
We denote by $(\pi,H,\Omega)$ the  GNS triple for $\varphi^H$. 
We set $M=\pi({\mathcal O}_\gamma(K))''$. 
Let $\hat{\varphi}^H$ be the normal extension of $\varphi^H$ to $M$ given by  
$\hat{\varphi}^H(m) = (m\Omega|\Omega)$ for $m \in M$. 
Since $\varphi^H$ is a KMS state, $\hat{\varphi}^H$ is faithful on $M$. 

Let $\tau ^H$ be the trace on $A=C(K_N)$ given by $\mu^H$ and  
let $(\pi_0,H_0,\Omega _0)$ be the  GNS triple for $\tau^H$. 
We can identify $H_0$ with $L^2(K,\mu ^H)$ 
and $\pi_0(a)$ with the multiplication operator by $a \in A$. 
Thus $\pi_0$ extends to a normal representation of $L^\infty(K,\mu^H)$. 
On the other hand, since $\tau^H$ is the restriction of $\varphi^H$ to $A$, 
the Hilbert space $H_0$ is naturally identified with the closure of 
$\pi(A)\Omega$ and the restriction $\pi|_A$ of $\pi$ to $A$ is 
quasi-equivalent to $\pi_0$ \cite[Lemma 4.1]{I}. 
Therefore $\pi|_A$ extends to a normal representation 
of $B:=L^\infty(K,\mu^H)$, which is denoted by $\hat{\pi}$. 

For a measurable function $f$ on $K$, we denote by $||f||_p$ the $L^p$-norm 
of $f$ in $L^p(K,\mu^H)$. 

For $m \in M$, we set 
$||m||_2=||m\Omega_L||$. 
It is well-known that the strong operator topology on the 
unit ball $M^1$ of $M$ coincides with the topology given by $||\cdot||_2$, 
and $M^1$ is complete with respect to $||\cdot||_2$. 

Let $\{h_n\}_{n=1}^\infty$ be an increasing sequence of non-negative 
functions in $A$ such that $h_n(x)=0$ for $x\in \cC(\gamma)$ 
and $\{h_n(y)\}_{n=1}^\infty$ converges to 1 for all 
$y\in K\setminus \cC(\gamma)$. 
For  $1\leq j\leq N$, we introduce $u_{j,n}\in X$ by setting 
$u_{j,n}(\gamma_j(x),x)=h_n(x)$ and 
$u_{j,n}(\gamma_k(x),x)=0$ for $k\neq j$.

\begin{lemma} Let the notation be as above. 
Then the sequence $\{\pi(S_{u_{j,n}})\}_{n=1}^\infty$ converges to 
an isometry, say $\tS_j$, in the $*$ strong operator topology. 
The isometries $\{\tS_j\}_{j=1}^N$ satisfy the Cuntz algebra relation.  
For the modular automorphism group $\{\sigma^{\hat{\varphi}^H}_t\}_{t\in \R}$, 
the following holds:   
$\sigma^{\hat{\varphi}^H}_t(\tS_j)=e^{-t\log N\sqrt{-1}}\tS_j$, 
$\forall t\in \R.$
\label{lemma;IIIH}
\end{lemma}

\begin{proof} By straightforward computation, we have $||S_{u_{j,m}}||\leq 1$ and 
$$||\pi(S_{u_{j,m}}-S_{u_{j,n}})||_2^2=||h_m-h_n||_2^2,$$
$$||\pi(S_{u_{j,m}}^*-S_{u_{j,n}}^*)||_2^2=\frac{1}{N}||h_m-h_n||_2^2. $$
Since we assume that $\mu(C(\gamma))=0$, these imply that 
$\{\pi(S_{u_{j,n}})\}_{n=1}^\infty$ and 
$\{\pi(S_{u_{j,n}})^*\}_{n=1}^\infty$ are 
Cauchy sequences with respect to 
$||\cdot||_2$, and so $\{\pi(S_{u_{j,n}})\}_{n=1}^\infty$ converges in the 
$*$ strong operator topology. 
As we have $S_{u_{j,n}}^*S_{u_{k,n}}=\delta_{j,k}\pi(|h_n|^2)$ which converges 
to $\delta_{j,k}1$ in the strong operator topology, the operators 
$\{\tS_j\}_{j=1}^N$ are isometries with mutual orthogonal ranges. 
Since $\sigma^{\hat{\varphi}^H}_t$ is a normal extension of 
$\alpha_{-t\log N}$, we have 
$\sigma_{t}^{\hat{\varphi}^H}(\tS_j)=e^{-t\log N\sqrt{-1}}\tS_j$, 
which implies 
$$\sum_{j=1}^N\hat{\varphi}^H(\tS_j\tS_j^*)
=\sum_{j=1}^N\frac{1}{N}\hat{\varphi}^H(\tS_j^*\tS_j)=1.$$
Since $\hat{\varphi}^H$ is faithful, this implies 
$\sum_{j=1}^N\tS_j\tS_j^*=1.$
\end{proof}

\begin{thm} Let the notation be as above. 
Then $M$ is the AFD type III$_{1/N}$ factor. 
\label{theorem;IIIH}
\end{thm}

\begin{proof} Thanks to the above lemma, it suffices to show that 
$M$ is generated by $\{\tS_j\}_{j=1}^N$.  

Since $aS_{u_{j,n}}=S_{u_{j,n}}a\cdot\gamma_j$ holds for all $a\in A$, we get 
$\pi(a)\tS_j=\tS_j\pi(a\cdot\gamma_j)$.   
For $f\in X$, we define $f_j\in A$ by 
$f_j(x)=f(\gamma_j(x),x)$. 
Since $\mu(C(\gamma))=0$, we have $\tS_j^*\pi(S_f)=\pi(f_j)$, and so
$$\pi(S_f)=\sum_{j=1}^N\tS_j\pi(f_j).$$ 
Therefore the linear span of the elements of the form 
$\tS_\eta \pi(a) \tS_\zeta^*$ with $a\in C(K)$, 
$\eta\in \{1,2,\cdots,N\}^m$, and $\zeta\in \{1,2,\cdots,N\}^n$ is a dense 
$*$-algebra of $M$, 
where $\tS_\eta=\tS_{\eta_1}\tS_{\eta_2}\cdots \tS_{\eta_m}$.  
We understand $\tS_\eta f\tS_\zeta^*=f$ for $\eta=\zeta=\emptyset$. 
This observation shows that to finishes the proof, it suffices to prove 
$\pi(A)\subset {\{\tS_j\}_{j=1}^N}''$.

For $a\in A$, we set 
$$a_n=\sum_{\eta\in \{1,2,\cdots,N\}^n}\hat{\varphi}^H(\tS_\eta^*\pi(a)
\tS_\eta)\tS_\eta\tS_\eta^*
=\sum_{\eta\in \{1,2,\cdots,N\}^n}\hat{\varphi}^H(\pi(a\cdot\gamma_\eta))
\tS_\eta\tS_\eta^*.$$
Then $\{a_n\}_{n=1}^\infty$ is a bounded sequence in ${\{\tS_j\}_{j=1}^N}''$. 
We show that $\{a_n\}_{n=1}^\infty$ converges to $\pi(a)$ in the strong 
operator topology. 
Indeed, we have 
\begin{eqnarray*}
||\pi(a)-a_n||_2^2&=&\hat{\varphi}^H(\pi(|a|^2)
-\pi(a^*)a_n-a_n^*\pi(a)+|a_n|^2)\\
&=&||a||_2^2-\frac{1}{N^n}
\sum_{\eta\in \{1,2,\cdots,N\}^n}|(f\cdot\gamma_\eta|1)_{\nu_N}|^2. 
\end{eqnarray*} 
Lemma \ref{lemma;IIIH} shows that the right-hand side tends to 0 as 
$n$ goes to infinity. 
Thus $\pi(A)\subset {\{\tS_j\}_{j=1}^N}''$. 
\end{proof}

\subsection{The case of Lyubich measures}
Let $R$ be a rational function with $N = \deg R \geq 2$, 
$A=C(\hat{\C})$, and $X=C(\graph R)$ be as before. 
Since $\hat{\C}\ni z\mapsto (z,R(z))\in \graph R$ is a homeomorphism, 
we identify $X$ with $C(\hat{\C})$, whose $A-A$ bimodule structure is given 
by 
$$(a\cdot f\cdot b)(z)=a(z)f(z)b(R(z)),$$
$$(f|g)_A(z)=\sum_{w\in R^{-1}(z)}e(w)\overline{f(w)}g(w).$$

The following theorem relies on Heicklen-Hoffman's remarkable 
result \cite{HH} saying that $(\hat{\C},R,\mu^H)$ is conjugate to 
$(K_N,\sigma,\nu_N)$ as measurable dynamical systems. 

\begin{thm} Let $R$ be a rational function with 
$N = \deg R \geq 2$. Let $\mu^L$ be the Lyubich measure 
and $\varphi^L$ be the corresponding 
$\log N$-KMS state on ${\mathcal O}_R(\hat{\mathbb C})$ 
for the gauge action. Then the GNS representation of $\varphi^L$ is 
a factor representation of type $III_{1/N}$. 
\label{theorem;IIIL}
\end{thm}

\begin{proof}
Let $(\pi,H,\Omega)$ be the  GNS triple for $\varphi^L$. 
We set $M=\pi({\mathcal O}_R(\hat{\mathbb C}))''$. 
Let $\hat{\varphi}^L$ be the normal extension of $\varphi^L$ to $M$. 
We denote by $\hat{\pi}$ the normal extension of the restriction $\pi|_A$ 
to $B:=L^\infty(\hat{\C},\mu^L)$. 
For $m \in M$, we set $||m||_2=||m\Omega_L||$. 
We denote by $||\cdot||_p$ the $L^p$-norm of 
$L^p(\hat{\C},\mu^L)=L^p(J_R,\mu^L)$. 

For $f\in X$, the condition  $\mu^H(\cC(R))=0$ implies  
$$||f||_{\infty} \leq ||\pi_L(S_f)||=||\pi_L((f|f)_A)||^{1/2}
\leq \sqrt{N}||f||_\infty,$$
$$||\pi_L(S_f)||_2^2  =
\int_{\hat{\C}}\sum_{w\in R^{-1}(z)}e(w)|f(w)|^2d\mu^L(z)
=\int_{\hat{C}} \sum_{w \in R^{-1}(z)}|f(z)|^2 d\mu ^L(z)
=N||f||_2^2.$$
The KMS condition implies $||\pi(S_f^*)||_2=||f||_2/\sqrt{N}$. 

For $f\in L^\infty(\hat{\C},\mu^L)$, we choose a sequence 
$\{f_n\}_{n=1}^\infty$ in $C(\hat{\C})$ satisfying the following 
three conditions: (1) $\{||f_n||_\infty\}_{n=1}^\infty$ is dominated by 
$||f||_\infty$, (2) $\{||f_n-f||_2\}_{n=1}^\infty$ converges to 0, and 
(3) $\{f_n(z)\}_{n=1}^\infty$ 
converges to $f(z)$ for almost every $z$ with respect to $\mu^L$. 
Then $\{\pi_L(S_{f_n})\}_{n=1}^\infty$ converges in the $*$ strong operator 
topology, whose limit is denoted by $\hat{S}_f$. 
Note that $\hat{S}_f$ does not depends on the choice of the sequence 
$\{f_n\}_{n=1}^\infty$. 
Let $g\in L^\infty(\hat{\C},\mu^L)$ and $\{g_n\}_{n=1}^\infty$ be a sequence 
in $C(\hat{\C})$ satisfying the same properties as above for $g$ instead of 
$f$.
Then 
$$\hat{S}_f^*\hat{S}_g=\lim_{n\to\infty}\pi_L((f_n|g_n)_A),$$
where the limit is taken in the $*$ strong operator topology. 
For $f,g\in B$, we define $(f|g)_B\in B$ by the same formula as the 
$A$-valued inner product of $X$. 
Then we get 
\begin{equation}
\hat{S}_f^*\hat{S}_g=\hat{\pi}((f|g)_B),\quad f,g\in 
L^\infty(\hat{\C},\mu^L).
\label{L1}
\end{equation}
In a similar way, we have 
\begin{equation}
\hat{\pi}(a)\hat{S}_f\hat{\pi}(b)=
\hat{S}_{afb\cdot R},\quad f,a,b\in 
L^\infty(\hat{\C},\mu^L).
\label{L2}
\end{equation}
Since $\sigma_t^{\hat{\varphi}^L}$ is a normal 
extension of $\alpha_{-t\log N}$, we have  
\begin{equation}
\sigma_{t}^{\hat{\varphi}^L}(\hat{\pi}(a))=\hat{\pi}(a),
\quad a\in L^\infty(\hat{\C},\mu^L). 
\label{L3}
\end{equation} 
\begin{equation}
\sigma_{t}^{\hat{\varphi}^L}(\hat{S}_f)=e^{-t\log N\sqrt{-1}}
\hat{S}_f,\quad f\in L^\infty(\hat{\C},\mu^L). 
\label{L4}
\end{equation}

Thanks to Heicklen and Hoffman \cite{HH}, there exists a 
von Neumann algebra isomorphism $\phi:L^\infty(K_N,\nu_N)\rightarrow 
L^\infty(\hat{\C},\mu^L)$ satisfying 
$\phi(f)\cdot R=\phi(f\cdot \sigma)$ for all 
$f\in L^\infty(K_N,\nu_N)$. 
We claim that there exists a representation $(\tilde{\pi},H)$ 
of $\cO_{\rho}(K_N)$ satisfying 
$$\tilde{\pi}(a)=\hat{\pi}(\phi(a)),\quad a\in C(K_N).$$
$$\tilde{\pi}(S_f)=\hat{S}_{\phi(f)},\quad f\in Y=C(K_N).$$
Thanks to (\ref{L1}) and (\ref{L2}), such a representation 
would exist for the Cuntz-Toeplitz algebra $\cT_Y$. 
Since $Y$ has a finite basis $\{1_{E_j}\}_{j=1}^N$, to prove 
the claim it suffices to show 
$$\sum_{j=1}^N\hat{S}_{\phi(1_{E_j})}\hat{S}_{\phi(1_{E_j})}^*=I.$$
Indeed, thanks to the KMS condition we have 
$$\hat{\varphi^H}(I-
\sum_{j=1}^N\hat{S}_{\phi(1_{E_j})}\hat{S}_{\phi(1_{E_j})}^*)
=1-\frac{1}{N}\sum_{j=1}^N
\hat{\varphi}^H(\hat{S}_{\phi(1_{E_j})}^*\hat{S}_{\phi(1_{E_j})})=0.$$ 
Since $\hat{\varphi}^H$ is faithful, we get the claim. 

It is routine work to show that $\tilde{\pi}(\cO_\rho(K_N))''=M$, 
and so $(\tilde{\pi},H, \Omega)$ is a cyclic representation. 
(\ref{L3}) and (\ref{L4})  imply that 
$(\tilde{\pi},H,\Omega)$ is unitarily equivalent  
to the GNS representation of the unique $\log N$-KMS state 
of the gauge action on $\cO_\rho(K_N)$. 
Thus Theorem~\ref{theorem;IIIH} shows that 
$M$ is  the AFD type III$_{1/N}$ factor.   
\end{proof}


\section{Symmetry} 
Since our construction of the Cuntz-Pimsner algebra from a dynamical system 
is natural, a symmetry of the dynamical system gives rise to an 
automorphism of the algebra. 
In this section, we give a criterion for outerness of such automorphisms. 
Quasi-free automorphisms on the Cuntz-Pimsner algebras for bimodules with 
finite bases are studied by Katayama-Takehana \cite{KT} and 
Zacharias \cite{Z} while the bimodules we treat in this paper do not 
necessarily have finite bases.  

Let $M$ be the AFD type III$_{1/N}$ factor acting on a Hilbert space. 
Since $M$ arises from the GNS representation of the unique KMS state 
of the gauge action on the Cuntz algebra $\cO_N$, 
we may find a generating set $\{S_j\}_{j=1}^N$ of $M$ consisting of isometries 
satisfying the following conditions: 
(1) they satisfy the Cuntz algebra relation and (2) there exists a cyclic and 
separating vector $\Omega\in H$ for $M$ such that if we denote 
the vector state for $\Omega$ by $\varphi$, then 
$\sigma^\varphi_t(S_j)=N^{-t\sqrt{-1}}S_j$ holds for $t\in \R$ and 
$j=1,2,\cdots,N$. 
The canonical shift endomorphism $\Phi$ of $M$ is defined by 
$$\Phi(x)=\sum_{j=1}^NS_jxS_j^*, \quad x\in M.$$
Let $B$ be the maximal abelian subalgebra of $M$ generated by 
$\bigcup_{n=1}^\infty \{S_\xi S_{\xi}^*\}_{\xi\in \{1,2,\cdots,N\}^n}.$

\begin{lemma} Let the notation be as above. 
Assume that $\theta$ is an automorphism of $M$ commuting with the modular 
automorphism group $\{\sigma^\varphi_t\}_{t\in \R}$ such that 
$\theta(B)=B$ and $\Phi\cdot\theta(x)=\theta\cdot \Phi(x)$ for all $x\in B$. 
If there exists $b\in B$ satisfying $\theta(b)\neq b$, then $\theta$ is outer. 
\label{lemma;outer}
\end{lemma}

\begin{proof} We first fix the notation. 
Let $\cF_n=\Phi^n(M)'\cap M$, which is isomorphic to 
$M_{N^n}(\C)$. 
The centralize $M_\varphi$ is a type II$_1$ factor generated by 
$\bigcup_{n=1}^\infty\cF_n$ and $M_\varphi$ satisfies $M\cap M_\varphi'=\C$. 
$B$ is a maximal abelian subalgebra of $M_\varphi$ as well. 

Suppose that there exists a unitary $u\in M$ satisfying 
$\theta(x)=uxu^*$ for all $x\in M$. 
We claim $u\in M_\varphi$. 
For $n\in \Z$, we set 
$$u_n=\frac{1}{T}\int_{0}^{T}N^{tn\sqrt{-1}}\sigma^\varphi_t(x)dt,$$
where $T=2\pi/\log N$. 
Then $u_mx=\theta(x)u_m$ holds for all $x\in M_\varphi$, which shows 
$u_m^*u_m, u_mu_m^*\in M\cap M_\varphi'=\C$. 
Thus $u_m$ is a multiple of a unitary. 
If $u_m\neq 0$ for some $m\neq 0$, the KMS condition would imply 
$\varphi(u_m^*u_m)\neq \varphi(u_mu_m^*)$, which is a contradiction. 
Thus we get $u=u_0\in M_\varphi$. 

Let $E_n$ be the trace preserving conditional expectation from 
$M_\varphi$ onto $\cF_n$. 
Then $\{||u-E_n(u)||_2\}_{n=1}^\infty$ converges to 0. 
For every $b\in B$, we have 
\begin{eqnarray*}||b-\theta(b)||_2&=&||\Phi^n(b)-\Phi^n(\theta(b))||_2
=||\Phi^n(b)-\theta(\Phi^n(b))||_2\\
&=&||\Phi^n(b)u-u\Phi^n(b)||_2=||\Phi^n(b)(u-E_n(u))-(u-E_n(u))\Phi^n(b)||_2\\
&\leq& 2||u-E_n(u)||_2||b||,
\end{eqnarray*}
which shows $b=\theta(b)$. 
This is a contradiction and $\theta$ is outer. 
\end{proof}

Let $R$ be a rational function with $\deg R \geq 2$ and 
let $\theta$ be a homeomorphism of $\hat{\C}$ commuting with $R$. 
Then by naturality, we get an automorphism of $\cO_{R}(\hat{\C})$ and 
$\cO_R(J_R)$, which we denote by $\theta_{\hat{\C}}$ and 
$\theta_{J_R}$ respectively. 
Namely, we define the actions of $\theta$ on $a\in C(\hat{\C})$ and 
$f\in X=C(\graph R)$ by $\theta\cdot a(z)=a(\theta^{-1}(z))$ and by 
$\theta\cdot f(z,R(z))=f(\theta^{-1}(z),R(\theta^{-1}(z)))$ respectively. 

\begin{thm} Let the notation be as above and let $\mu^L$ be the 
Lyubich measure. 
If the automorphism of $L^\infty(\hat{\C},\mu^L)=L^\infty(J_R,\mu^L)$ 
induced by $\theta$ is non-trivial, then 
$\theta_{\hat{\C}}$ and $\theta_{J_R}$ are outer. 
\label{theorem;outer}
\end{thm}

\begin{proof} We use the notation in the proof of Theorem~\ref{theorem;IIIL}. 
Since $\theta_{\hat{\C}}$ preserves $\varphi^L$, there exists an 
automorphism $\theta_M$ of $M$ satisfying 
$\theta_M\cdot\pi=\pi\cdot \theta_{\hat{\C}}$. 
Since $\pi$ factors through $\cO_R(J_R)$, to prove the theorem 
it suffices to show that $\theta_M$ is outer. 

Let $S_j=\hat{S}_{\phi(1_{E_j})}$ and 
$$\Phi(x)=\sum_{j=1}^N S_jxS_j^*,\quad x\in M.$$
We claim that for any countable basis $\{u_n\}_{n=1}^\infty$ of $X$ and 
$f\in L^\infty(\hat{\C},\mu^L)$, the following hold:
$$\Phi(\hat{\pi}(f))=\sum_{j=1}^\infty\pi(S_{u_j})\hat{\pi}(f)\pi(S_{u_j}^*).$$
Indeed, thanks to Lemma~\ref{lemma;p=1}, the right-hand side converges in 
the strong operator topology. 
Since $S_k^*\pi(S_{u_j})$ commutes with $\hat{\pi}(A)''$, 
we have 
\begin{eqnarray*}\pi(S_{u_j})\hat{\pi}(f)(S_{u_j}^*)&=&\sum_{k=1}^N
S_kS_k^*\pi(S_{u_j})\hat{\pi}(f)\pi(S_{u_j}^*)
=\sum_{k=1}^N
S_k\hat{\pi}(f)S_k^*\pi(S_{u_j})\pi(S_{u_j}^*)\\
&=&\Phi(\hat{\pi}(f))\pi(S_{u_j})\pi(S_{u_j}^*),
\end{eqnarray*}
which shows the claim. 

Since $\{\theta_{\hat{\C}}(u_j)\}_{j=1}^\infty$ is a basis of $X$ too, 
we get $\theta_M(\Phi(\hat{\pi}(f)))=\Phi(\theta_M(\hat{\pi}(f)))$ 
for every $f\in L^\infty(\hat{\C},\mu^L)=\phi(L^\infty(K_N,\nu_N))$. 
Now the statement follows from Lemma~\ref{lemma;outer} as 
the proof of Theorem~\ref{theorem;IIIH} shows that $\pi(A)''$ coincides 
with $B$ in Lemma~\ref{lemma;outer}. 
\end{proof}

\begin{exam}
We denote by $G_R$ the group of M\"{o}bius 
transformations commuting with $R$. 
Then $G_R$ is a finite subgroup of $PSL(2,{\mathbb C})$ 
(see  \cite[pp.\ 104]{B}) and $G_R$ is either a cyclic group, 
a dihedral group or the group of symmetries of 
a regular tetrahedron, octahedron or icosahedron.
It is known that 
there exist exactly three rational functions 
with $2 \leq \deg R \leq 30$ such that 
$G_R$ is the group of symmeties of a regular 
icosahedron, which is isomorphic to $A_5$ 
(see Doyle and McMullen \cite{DMc}). 

Let $R(z) = z^n$ for $n \geq 2$ and let $\omega$ be the primitive $(n-1)$th 
root of the unity. Then 
\[
G_R = \{z, \omega z, \dots, \omega ^{n-2}z, 
\frac{1}{z}, \frac{\omega}{z}, \dots, \frac{\omega^{n-2}}{z} \}. 
\]
Thus $G_R$ is isomorphic to the dihedral group $D_{n-1}$ for 
$n \geq 3$  and is isomorphic to ${\mathbb Z}/2{\mathbb Z}$ for 
$n = 2$.  
Therefore the $G_R$-action on $\cO_{R}(\hat{\C})$ transitively acts on 
the set of extreme $\beta$-KMS states for $\beta > \log n$. 
\end{exam}

\begin{exam} If $R$ has only real coefficients, 
then $\theta(z)=\overline{z}$ commutes with $R$. 
\end{exam}

\begin{rema} Let $\gamma=(\gamma_1,\gamma_2,\cdots,\gamma_N)$ be a system 
of proper contractions on a compact metric space $K$ satisfying 
the self-similarity condition.  
Then for a homeomorphism $\theta$ of $K$ satisfying 
$\theta\cdot\gamma_j\cdot\theta^{-1}=\gamma_{p(j)}$, where 
$p$ is a permutation of $\{1,2,\cdots,N\}$, a result 
similar to Theorem~\ref{theorem;outer} follows from Lemma~\ref{lemma;outer}.  
\end{rema}

\end{document}